\documentstyle[amssymb,12pt]{article}

\input amssym.def
\input amssym

\newtheorem{thm}{Theorem}[section]
\newtheorem{lem}[thm]{Lemma}

\newtheorem{cor}[thm]{Corollary}
\newtheorem{fact}[thm]{Fact}

\newtheorem{defn}[thm]{Definition}

\newtheorem{nrmk}[thm]{Remark}

\newcommand{\qed}{\hfill $\Box$ \vspace{.5cm}}
\newcommand{\pf}{{\bf Proof. }}

\title {Solvable groups definable in o-minimal structures}

\author {M\'{a}rio J. Edmundo 
\\The Mathematical Institute
\\24-29 St Giles
\\OX1 3LB  Oxford, U.K
\\edmundo@maths.ox.ac.uk\\}
\date{January 24, 2000}

\newcommand{\into}{\longrightarrow}

\renewcommand{\hat}{\widehat}
\renewcommand{\tilde}{\widetilde}
\renewcommand{\bar}{\overline}

\newcommand{\NN}{\mathbb{N}}

\newcommand{\RR}{\mathbb{R}}

\newcommand{\M}{\mbox{$\cal M$}}
\newcommand{\N}{\mbox{$\cal N$}}

\newcommand{\I}{\mbox{$\cal I$}}
\newcommand{\J}{\mbox{$\cal J$}}

\begin{document}

\maketitle
\begin{abstract}
Let $\N$ be an o-minimal structure. In this paper  
we develop group extension theory 
over $\N$ and use it to describe $\N$-definable
solvable groups. We prove an o-minimal analogue of the
Lie-Kolchin-Mal'cev theorem and we describe $\N$-definable $G$-modules
and $\N$-definable rings.
\end{abstract}

\newpage

\begin{section}{Introduction}
\label{section introduction}

We will work inside an o-minimal structure $\N$$=(N,<,\ldots )$ and 
therefore definable will mean $\N$-definable. We will assume the readers
familiarity with basic o-minimality (see \cite{vdd}). We will start by
recalling some basic notions and results on definable groups that will
be used through the paper.

\medskip
Pillay in \cite{p} adapts Hrushovski's proof of Weil's Theorem that an
algebraic group can be recovered from birational data to show 
that a definable group $G$ can be equipped
with a unique definable manifold structure making the group into a topological
group, and definable homomorphisms between 
definable groups are topological homomorphisms. 
In fact, as remarked in \cite{pps1}, if $\N$ is an o-minimal
expansion of a real closed field $G$ equipped with the above unique 
definable manifold structure is a $C^p$ group for all $p\in \NN$; 
and definable homomorphisms between 
definable groups are $C^p$ homomorphisms for all $p\in \NN$ .
Moreover, again by \cite{pps1}, the definable manifold structure on a
definable subgroup is the sub-manifold structure.

By \cite{p} definable groups satisfies the descending chain condition 
(DCC) on definable subgroups. This is used to show that 
the definably connected component of the identity $G^0$ of a definable group
$G$ is the smallest definable subgroup of $G$ of finite index. Also
infinite such groups have  infinite definable abelian subgroups; a
definable subgroup $H$ of $G$ is closed and the following are
equivalent (i) $H$ has finite index in $G$, (ii) $dimH=dimG$, (iii)
$H$ contains an open neighbourhood of the identity element of $G$ and (iv)
$H$ is open in $G$. Finally, by \cite{s} an infinite abelian definable
group $G$ has unbounded
exponent and the subgroup $Tor(G)$ of torsion points of $G$ is
countable. In particular, if $\N$ is $\aleph _0$-saturated then $G$ has
an element of infinite order.

\medskip
One dimensional definable manifolds are classified in \cite{r}
and the following is deduced. Suppose that $G$ is one-dimensional
definably connected definable group. Then  by \cite{p} $G$ is abelian, and 
$G$ is torsion-free or for each prime $p$ the set of $p$-torsion
points of
$G$ has $p$ elements. In the former case $G$ is an ordered abelian
divisible definably simple group.

Note that if $I$ is a one-dimensional definably connected 
ordered definable group,
then the structure $\I$ induced by $\N$ on $I$ is o-minimal.
In particular, we have the following results 
from \cite{ms}. Suppose that $(I,0,1,+,<)$ is a one-dimensional
definably connected torsion-free definable group, where $1$ is a 
fixed positive element. Let $\Lambda(\I)$ be
the division ring of all $\I$-definable endomorphisms of $(I,0,+).$
Then exactly one of the following holds: (1) $\I$ is {\it linearly bounded}
with respect to $+$ (i.e, for every
$\I$-definable function $f:I\into I$ there is $r\in \Lambda (\I)$ such
that $\lim _{x\into +\infty}[f(x)-rx]\in I$), or (2) there is a  
$\I$-definable binary
operation $\cdot $ such that $(I,0,1,+,\cdot,<)$ is a real closed
field. Also, up to $\I$-definable isomorphism there is at most one
$\I$-definable group $(I,0,* )$ such that $\I$ is linearly
bounded with respect to $* $ and at most one $\I$-definable
(real closed) field $(I,0,1,\oplus ,\otimes )$.

Moreover, the following are equivalent: (i)
$\I$ is linearly bounded with respect to $+$, (ii) for every
$\I$-definable function $f:A\times I\into I$, where $A\subseteq I^n$,
there are $r_1,\dots ,r_l\in \Lambda(\I)$ such that for every $a\in A$
there is $i\in \{1,\dots ,l\}$ with 
$\lim _{x\into +\infty}[f(a,x)-r_ix]\in I$ and (iii) there is no
infinite definable subset of $\Lambda (\I)$.

Let $(I,0,1,+,<)$ be as above and let $\Lambda :=\Lambda (\I)$. 
Then $\I$ is called {\it semi-bounded}
if every $\I$-definable set is already definable in the reduct
$(I,0,1,+,<,(B_k)_{k\in K},(\lambda )_{\lambda \in \Lambda }$$),$
of $\I$ where $(B_k)_{k\in K}$ is the collection of all bounded
$\I$-definable sets. According to \cite{e}, the following are
equivalent: (i) $\I$ is semi-bounded, (ii) there is no $\I$-definable
function between a bounded and an unbounded subinterval of $I$, (iii)
there is no $\I$-definable (real closed) field with domain an
unbounded subinterval of $I$, (iv) for every $\I$-definable function
$f:I\into I$ there are $r\in \Lambda $, $x_0\in I$ and $c\in I$
such that for all $x>x_0$, $f(x)=rx+c$ and (v) $\I$ satisfies the ``structure
theorem''.

Let $(I,0,1,+,\cdot ,<)$ be a real closed field definable in $\N$. Let
$\mathcal{K}(\I)$ be the ordered field of all
$\I$-definable endomorphisms of the multiplicative group
$(I^{>0},\cdot ,1)$. Note that $\mathcal{K}(\I)$$\into I$, $\alpha
\into \alpha '(1)$ is an embedding of ordered fields.
The elements of $\mathcal{K}(\I)$ are called {\it power functions}
and for $\alpha \in \mathcal{K}(\I)$ with $\alpha '(1)=r$ we write
$\alpha (x)=x^r$.
By \cite{m} exactly one of the following holds: (1) $\I$ is {\it
power bounded} (i.e., for every $\I$-definable function $f:I\into I$
there is $r\in \mathcal{K}(\I)$ such that ultimately $|f(x)|<x^r$) or (2) 
$\I$ is {\it exponential} (i.e., there is an $\I$-definable ordered group
isomorphism $e:(I,0,+,<)\into (I^{>0},1,\cdot ,<)$). Moreover, the
following are equivalent: (i) $\I$ is power bounded, (ii) for every
$\I$-definable function $f:A\times I\into I$, where $A\subseteq I^n$,
there are $r_1,\dots ,r_l\in \mathcal{K}(\I)$ such that for every $a\in
A$, if the function $x\into f(a,x)$ is ultimately nonzero then, 
there is $i\in \{1,\dots ,l\}$ with 
$\lim _{x\into +\infty}[f(a,x)/x^{r_i}]\in I$ and (iii) there is no
infinite definable subset of $\mathcal{K}(\I)$.

If $\I$ is power bounded, then we know that $(I,0,+,<)$ and
$(I^{>0},1,\cdot ,<)$ are the only (up to $\I$-definable isomorphism) 
$\I$-definable one-dimensional
torsion-free ordered groups. The {\it Miller-Starchenko conjecture}
says that in an o-minimal expansion $\I$ of an
ordered field every $\I$-definable one-dimensional
torsion-free ordered group is $\I$-definable isomorphic to either $(I,0,+,<)$
or $(I^{>0},1,\cdot ,<)$. (In the general case we only know (see
\cite{ms}) that up to $\I$-definable isomorphisms there are at most
two $\I$-definably connected, $\I$-definable one-dimensional
torsion-free ordered groups). 
Suppose that the Miller-Starchenko conjecture does not
hold for $\I$, then we call the unique $\I$-definable group $(I,0,\oplus
,<)$ which is not $\I$-definably isomorphic to $(I,0,+,<)$ or
$(I^{>0},1,\cdot ,<)$ the {\it Miller-Starchenko group of $\I$}.
Note the following: if $G$ is an $\I$-definable
one-dimensional torsion-free ordered group, then we can assume that
$G=(I,0,\oplus ,<)$, and $\alpha :G\into (I,0,+)$ is an abstract $C^1$ 
isomorphism iff $\forall
s\in G,\,\alpha '(s)\frac{\partial \oplus}{\partial x}(0,s)=\alpha '(0)$ 
where for all $t,s\in G,\, \oplus (t,s):=t\oplus s$ i.e., $\alpha $ is
Pfaffian over $(I,0,1,+,\cdot,\oplus ,<)$ (note that, by associativity
of $\oplus $, for all $s\in G$, $\frac{\partial \oplus}{\partial
x}(0,s)\neq 0$). 

\medskip
The notion of definably compact groups was introduced in \cite{ps}.
Let $G$ be a definable group. We say that $G$ is {\em definably
compact} if for every definable continuous embedding $\sigma\colon
(a,b)\subseteq N\to G$, where $-\infty\leq a<b\leq +\infty$, there are
$c,d\in G$ such that $\lim_{x\to a^+}\sigma(x)=c$ and
$\lim_{x\to b^-}\sigma(x)=d$, where the limits are taken with respect to the
topology on $G$. In \cite{ps} the following result is established. 
Let $G$ be a definable group which is not definably compact. Then $G$
has a one-dimensional definably connected torsion-free (ordered)
definable subgroup.

\medskip
The trichotomy theorem \cite{pst1} and the theory of non orthogonality
from \cite{pps1} are used to prove the following (see theorem \ref{thm
definable quotients} and theorem \ref{thm existence of definable sections}).

\begin{fact}\label{fact definable quotients}
Let $U$ be a definable group and let $A$ be a definable normal
subgroup of $U$. Then we have a definable extension
$1\rightarrow A\rightarrow U\stackrel{j}\rightarrow G\rightarrow 1$
with a definable section $s:G\into U$.
\end{fact}

If we take in fact \ref{fact definable quotients} $A$ to be the
definable radical of $U$ i.e., the maximal definable solvable normal 
subgroup of
$U$ we get that $G$ is either finite or definably semisimple i.e., it
has no infinite proper abelian definable normal subgroup. Definable
definably semisimple groups are classified in \cite{pps1} (see also
\cite{pps2} and \cite{pps3}). 
Below, $\mathbf{G}$ is the structure $(G,\cdot )$ where
$\cdot $ is the group operation of $G$.

\begin{fact}\label{fact semisimple}
\cite{pps1} and \cite{pps3}.
Let $G$ be a definably semisimple $\mathbf{G}$-definably connected
definable group. Then $G=G_1\times \cdots \times G_l$ and for each 
$i\in \{1,\dots ,l\}$ there is an o-minimal expansion $\I$$_i$ of a
real closed field definable in $\N$ such that there is no definable 
bijection between a distinct pair among the $I_i$'s,
$G_i$ is $\I$$_i$-definably isomorphic to a $I_i$-semialgebraic
subgroup of $GL(n_i,I_i)$ which is a direct product of
$I_i$-semialgebraically simple, $I_i$-semialgebraic subgroups of
$GL(n_i,I_i)$.
\end{fact} 

Fact \ref{fact definable quotients} allows us to develop group
extension theory with  abelian  and non abelian
kernel over $\N$. We use this theory to prove the following 
result for definable solvable groups (see theorem \ref{thm the main
theorem}).
 
\begin{fact}\label{fact the main theorem} 
Let $U$ be a definable solvable group. Then $U$ has a definable normal
subgroup $V$ such that $U/V$ is a definably compact definable solvable
group and 
$V=K\times W_1\times \cdots \times W_s
\times V'_1\times V_1 \times \cdots \times V'_k\times V_k$ 
where $K$ is the definably connected
definably compact normal subgroup of $U$ of maximal dimension and 
for each $j\in \{1,\dots ,s\}$ (resp., $i\in \{1,\dots ,k\}$) 
there is a semi-bounded o-minimal expansion $\J$$_j$ of a group 
(resp., an o-minimal expansion $\I$$_j$ of a real closed field)
definable in $\N$ such that there is no definable bijection between
a distinct pair among the $J_j$'s and $I_i$'s,
$W_j$ is a direct product of copies of the additive group of $\J$$_j$, 
$V'_i$ is a direct product of copies of the linearly bounded 
one-dimensional torsion-free $\I$$_i$-definable group and 
$V_i$ is an $\I$$_i$-definable group such that
$Z(V_i)$ has an $\I$$_i$-definable subgroup $Z_i$ such that
$Z(V_i)/Z_i$ is a direct product of copies of the linearly bounded
one-dimensional torsion-free $\I$$_i$-definable group and there are
$\I$$_i$-definable subgroups $1=Z^0_i<Z^1_i<\cdots <Z^{m_i}_i=Z_i$
such that for each $l\in \{1,\dots ,m_i\}$, $Z^l_i/Z^{l-1}_i$ is the
additive group of $\I$$_i$, and $V_i/Z(V_i)$ $\I$$_i$-definably embeds
into $GL(n_i,I_i)$.
\end{fact}

We also prove the following result about definably compact definable
groups (see corollary \ref{cor compact U}).

\begin{fact}\label{fact compact U}
Let $U$ be a definably compact, definably connected
definable group. Then $U$ is either abelian or
$U/Z(U)$ is a definable semi-simple group. 
In particular, if $U$ is solvable then it is abelian.
\end{fact}

Fact \ref{fact the main theorem} gives a partial solution to the
Peterzil-Steinhorn splitting problem for solvable definable group with
no definably compact parts (see \cite{ps}).
We say that a definable abelian group $U$ {\it has no definable 
compact parts} if there are definable subgroups
$1=U_0<U_1<\cdots <U_{n}=U$ such that for each 
$j\in \{1,\dots ,n\}$, 
$U_j/U_{j-1}$ is a one-dimensional definably connected
torsion-free definable group.
We say that a definable solvable group $U$ {\it has no definable 
compact parts} if $U$
has definable subgroups $1=U_0\trianglelefteq U_1\trianglelefteq
\cdots \trianglelefteq U_n=U$ such that for each $i\in \{ 1, \dots
,n\}$, $U_i/U_{i-1}$ is a definable abelian group with no definable
compact parts. Peterzil and Steinhorn ask in \cite{ps} if a definable
abelian group $U$ of dimension two and with no definably compact parts is
a direct product of one-dimensional definably connected
torsion-free definable groups. Fact  \ref{fact the main theorem} above
reduces this problem to the case where $U$ is a group definable in a
definable o-minimal expansion $\I$ of a real closed field
$(I,0,1,+,\cdot ,<)$ and we have an $\I$-definable extension
$1\rightarrow A\rightarrow U\rightarrow G\rightarrow 1$ where
$A=(I,0,+,<)$ and $G=(I,0,*,<)$ is a one-dimensional torsion-free
$\I$-definable group. We prove (see lemma \ref{lem ps problem}) that
in this case, there is an $\I$-definable $2$-cocycle $c\in
Z^2_{\mathcal{I}}$$(G,A)$ for $U$ such that $U$ is $\I$-definably
isomorphic to $A\times G$ iff there is an $\I$-definable function
$\alpha :G\into A$ such that $\forall s\in G,\, \alpha
'(s)\frac{\partial *}{\partial x}(0,s)=\alpha '(0)+\frac{\partial
c}{\partial x}(0,s)$.

Let $\I$ be an o-minimal expansion of a real closed field
$(I,0,1,+,\cdot ,<)$ and suppose that we have an abelian 
$\I$-definable extension
$1\rightarrow A\rightarrow U\rightarrow G\rightarrow 1$ where
$A=(I,0,+,<)$ and $G=(I,0,*,<)$ is a one-dimensional torsion-free
$\I$-definable group. We shall say that $U$ is a {\it
Peterzil-Steinhorn $\I$-definable group} if $U$ is not $\I$-definably
isomorphic to $A\times G$.
A corollary of our main result is the following fact (see corollary
\ref{cor minimal reduct}).

\begin{fact}\label{fact minimal reduct}
Let $\I$$=(I,0,1,+,\cdot ,<,\dots )$ be an o-minimal expansion of a
real closed field with no Peterzil-Steinhorn $\I$-definable
groups. Then every $\I$-definable solvable group $U$ with no $\I$-definable 
compact parts is $\I$-definably
isomorphic to a group definable of the form 
$U'\times G_1\cdots G_k\cdot G_{k+1}\cdots G_l$ where $U'$ is a direct product
of copies of linearly bounded one-dimensional torsion-free
$\I$-definable groups, for $i=1,\dots ,k$, $G_i=(I,0,+)$
and for $i=k+1,\dots ,l$, $G_i=(I^{>0},1,\cdot )$. In
particular, $G:=G_1\cdots G_k\cdot G_{k+1}\cdots G_l$ $\I$-definably
embeds into some $GL(n,I)$ and $U$ is $\I$-definably isomorphic to a
group definable in one of the following  reducts $(I,0,1,+,\cdot
,\oplus )$,  $(I,0,1,+,\cdot ,\oplus ,e^t)$ or 
$(I,0,1,+,\cdot ,\oplus , t^{b_1},\dots
,$$t^{b_r})$ of $\I$ where $(I,0,\oplus )$ is the Miller-Starchenko
group of $\I$, $e^{t}$ is the $\I$-definable exponential map (if it
exists), and the $t^{b_j}$'s are $\I$-definable power
functions. Moreover, if $U$ is nilpotent then $U$ is $\I$-definably 
isomorphic to a
group definable in the reduct $(I,0,1,+,\cdot ,\oplus )$ of $\I$. 
\end{fact}
 
In section \ref{section lie kolchin malcev theorem}
we use our main result to classify definable $G$-modules (see theorem 
\ref{thm non trivial}), this is then used to prove the o-minimal version of the
Lie-Kolchin-Mal'cev theorem (see theorem \ref{thm lie kolchin malcev}).

Another application of fact \ref{fact the main theorem} is the 
following result (see theorem \ref{thm definable choice}):
Let $U$ be a definable group and let $\{T(x):x\in X\}$ be a definable
family of non empty definable subsets of $U$. Then there is a definable
function $t:X\into U$ such that for all $x,y\in X$ we have
$t(x)\in T(x)$ and if $T(x)=T(y)$ then $t(x)=t(y)$.
This result shows that the many of the theorems from
\cite{pst2} can be obtained without the assumption that $\N$ has
definable Skolem functions. We include here direct proofs (avoiding
the use of $\bigvee $-definability theory) of some of
these results, namely fact \ref{fact compact U} above, corollary 
\ref{cor defining a field} and corollary \ref{cor homo compact}. 

In section \ref{section definable rings} we apply the main theorem 
to describe definable rings (see theorem 
\ref{thm definable rings} and theorem \ref{thm compact rings}).

\end{section}

\begin{section}{Definable quotients}
\label{section definable quotients}

\begin{defn}\label{defn definable choice}
{\em
Let $S$ be a definable set and let 
$T:=\{T(x):x\in X\}$ be a definable family
of non empty definable subsets of $S$.
We say that $T$ has {\it definable choice} if there is a definable 
function $t:X\into S$
such that for all $x\in X$, $t(x)\in T(x)$. If in addition, $t$ is such
that for all $x,y\in X$, if $T(x)=T(y)$ then $t(x)=t(y)$, then we 
say that $T$ has {\it strong definable choice}.
The function $t$ is called 
a {\it (strong) definable choice for the family $T$}.
We say that  the definable set $S$ has {\it (strong) definable choice}
if every definable family $T$ of non empty definable subsets of $S$
has a (strong) definable choice.
}
\end{defn}

The following fact is easy to prove.

\begin{fact}\label{fact definable choice}
{\em
The following hold: (i) if $f:R\into S$ is a definable
map such that for all $s\in S$, $f^{-1}(s)$ is finite and $S$ has 
(strong)
definable choice then $R$ has (strong) definable choice; (ii) if
$g:S\into R$ is a surjective definable map and $S$ has (strong)
definable choice then $R$ has (strong) definable choice; (iii) if
$S:=S_1\times \cdots \times S_k$ is definable and each $S_i$ is
definable and has (strong) definable choice then $S$ has (strong)
definable choice. 
}
\end{fact}

For the prove of the next lemma we need to recall some definitions 
from \cite{pps1}: an open interval $I\subseteq N$ is {\it transitive}
if for all $x,y\in I$ there are definably homeomorphic subintervals
$I_x,I_y$ of $I$ containing $x$ and $y$ respectively; an open
rectangular box $I_1\times \cdots \times I_n$ is transitive if all the
intervals $I_k$ are transitive.

\begin{lem}\label{lem local definable choice}
A definable group $U$ has a definable 
neighbourhood $O$ of $1$ (the identity) with strong definable choice.
\end{lem}

\pf
Since it is sufficient to prove the lemma for an $\omega _1$-
saturated elementary extension of $\N$, we will assume that
$\N$ is $\omega _1$-saturated.

By lemma 1.28\cite{pps1}, there is a definable chart $(O',\phi )$ 
on $U$ at $1$ such that $\phi (O')$ is a transitive rectangular box, 
say $I_1\times \cdots \times I_n$. Let $\phi (1):=(a_1,\dots ,a_n)$.
Then by the trichotomy theorem \cite{pst1}, the definable
structure $\J$$_i$ induced by $\N$ on some open
subinterval $J_i$ of $I_i$ containing $a_i$ is either an o-minimal expansion of
of a real closed field or an o-minimal expansion of an ordered partial
group. Without loss of generality we may assume that $(J_i,a_i,+,-)$
is a definable ordered partial group with zero $a_i$ and
$J_i=(-e,e)$. By fact ref{fact 
definable choice} its enough to show that
$J'_i=(-\frac{e}{2},\frac{e}{2})$ has strong definable
choice. This is follows from the fact that there are definable
functions $l,r:J'_i\into J'_i$
and $m:J'_i\times J'_i\into J'_i$ such that for all $x,y\in J'_i$, we have 
$l(x)<x$, $x<r(x)$ and if $x<y$ then $x<m(x,y)<y$: take
$l(x):=x+|\frac{e-x}{2}|$; $r(x):=x-|\frac{e-x}{2}|$ and 
$l(x):=x+|\frac{y-x}{2}|$. 
\qed

Recall that, if we have a definable set $S$ and a definable
equivalence relation $E$ on $S$ then, we say that $S/E$ is 
definable if there is a definable map $\pi :S\into T$ such
that $\forall x,y\in S,\, xEy\iff \pi (x)=\pi (y)$. 
Note that this is the case, if the definable family 
$\{x/E:x\in S\}$ has a strong definable choice. 
If $S$ is a definable group, $E$ a definable normal subgroup and the
set $S/E$ is definable then, $S/E$ becomes in a natural way a
definable group.

\begin{thm}\label{thm definable quotients}
Let $U$ be a definable group and let $V$ be a definable normal
subgroup of $U$. Then $U/V$ is definable.
\end{thm}

\pf
Suppose that $U\subseteq N^m$ and for each $l\in \{1,\dots ,m\}$
let $\pi _l:N^m\into N^l$ be the projection onto the first $l$
coordinates and let $\pi ^l:N^m\into N$ be the projection onto the
$l$-th coordinate. 

The existence of a strong definable
choice $l:=(l_1,\dots ,l_m)$ for the family $\{xV:x\in U\}$
follows from the claim below. In fact the claim implies the existence
of $l$ on a large definable subset $U_m$ of $U$ (i.e.,
$dim(U\setminus U_m)<dimU$), but by lemma 2.4 
\cite{p}, there are $u_1,\dots ,u_n\in U$ such that 
$U=u_1U_m\cup \cdots \cup u_nU_m$ and so we can extend $l$ from $U_m$ to 
$U$.

\medskip
{\it Claim:}
For each $k\in \{1,\dots ,m\}$ there is a definable subset $U_k$
of $U$ such that (i) 
$dim(U\setminus U_k)<dimU$ and (ii)
if $x\in U_k$ and $y\in U$ is such that $xV=yV$ then
$y\in U_k$. Moreover, there are definable functions
$l_1,\dots ,l_k:U_k\into N$ such that for each 
$x\in U_k$ there is $z\in xV$ such that 
$\pi _k(z)=(l_1(x),\dots ,l_k(x))$ and for all $y\in U$ if $xV=yV$ then
$(l_1(x),\dots ,l_k(x))=(l_1(y),\dots ,l_k(y))$.

\medskip
{\it Proof of Claim:}
We will do this by induction on $k$. Suppose that $k=1$.
For the induction  let us introduce the following notation:
$U_0:=U$ and for each $x\in U_0$, let $V_0(x):=xV$.

We have a definable function
$\alpha _1:U_0\into N\cup \{+\infty \}$ given by, for each $x\in U_0$, 
$\alpha _1(x)=\sup \pi ^1(V_0(x))$. 
Note that, if $V_0(x)=V_0(y)$ then $\alpha _1(x)=\alpha _1(y)$.
Now if $x\in U_0$ is such that $\alpha _1(x)\in \pi ^1(V_0(x))$ then we
can take $l_1(x):=\alpha _1(x)$. 
Let $U'_0=U_0\setminus M_1$ where $M_1:=
\{x\in U_0:\alpha _1(x)\in \pi ^1(V_0(x))\}$
and suppose that $U'_0$ is non empty.
By o-minimality, the set $I_1$ of end points of $\alpha _1(U'_0)$
in $\alpha _1(U_0)$ is finite.
Suppose that $I_1$ is non empty and let $a\in I_1$. 
Consider the definable sub family $\{V_0(x):\alpha _1(x)=a\}$ of 
$\{V_0(x):x\in U_0\}$. Let $x_0\in U_0$ such that $\alpha _1(x_0)=a$ 
and 
define for all $x\in U_0$ such that
$V_0(x)=V_0(x_0)$, $l_1(x):=a_0$ where $a_0$ is some fixed
element of $\pi ^1(V_0(x_0))$. 
For each $x\in U_0$ such that $\alpha _1(x)=a$ let $\gamma _1(x):=\inf
\{z:a_0\leq z<a,\, (z,a)\subseteq \pi ^1(V_0(x))\}$. If
$V_0(x)=V_0(y)$ then $\gamma _1(x)=\gamma _1(y)$. For $x\in U_0$ with 
$\alpha _1(x)=a$ let $K_1(x):=\{z\in O:\alpha _1(zx)\in (\gamma
_1(x),a)\}$ where $O$ is the definable neighbourhood of $1$ in $U$ with 
strong definable choice (see lemma \ref{lem local definable choice}). 
This is a definable family of definable non-empty sets such that if 
$V_0(x)=V_0(y)$ then $K_1(x)=K_1(y)$. On $\{x\in U_0:\alpha _1(x)=a\}$
define $l_1(x):=\alpha _1(k_1(x)x)$ where $k_1(x)$ is a strong
definable choice for $K_1(x)$.  

If $X_1\cup M_1$ is large in 
$U_0$ then the claim is proved for $k=1$. Otherwise, we have
$dim(U_0\setminus (X_1\cup M_1))=dimU_0$.
Now let 
$J_1:=\alpha _1(U_0)\setminus I_1$. Suppose that $J_1$ is
non empty. Then $J_1$ is a finite union of open intervals.
Let $Y_1$ be the definable set of all $x\in U_0$ 
such that $\alpha _1(x)\in J_1$ and there is
(equivalently, for all) $y\in U_0$ such that
$V_0(y)= V_0(x)$ and $\alpha _1$ is continuous at $y$.
O-minimality implies
that $Y_1$ is large in $U_0\setminus (X_1\cup M_1)$
and so, $Y_1\cup X_1\cup M_1$ is large in $U_0$. 

Let $A_1$ be the definable subset of $Y_1$ of all $x\in Y_1$ such that
there is a definable open neighbourhood $B$ of $x$ in $U$, such that
$\alpha _1(B)\subseteq \{z\in J_1:\alpha _1(x)\leq z\}$. 
If $V_0(x)=V_0(y)$ and $x\in A_1$ then $y\in A_1$. Clearly,
by o-minimality, $\alpha _1(A_1)$ is finite and as before we can 
construct $l_1$ on $A_1$. 

Let $B_1:=Y_1\setminus A_1$ and suppose that
$B_1$ is non empty.
Then we have a definable family 
$\{T_1(x):x\in B_1\}$ of definable
subsets of $O$, the definable neighbourhood of $1$ in $U$ with 
strong definable choice (see lemma \ref{lem local definable choice}) 
given by
$T_1(x):=\{z\in O:\alpha _1(zx)\in S_1(x)\}$ where
$S_1(x):=\pi ^1(V_0(x))\cap \{z\in J_1:z< \alpha _1(x)\}$. 
By construction, for all $x\in B_1$, $S_1(x)$ is infinite  and if
$V_0(x)=V_0(y)$, then $y\in B_1$, $S_1(x)=S_1(y)$ and $T_1(x)=T_1(y)$. 
We now show that $T_1(x)$ is infinite for all $x\in B_1$: let 
$z'<\alpha _1(x)$ such that $(z',\alpha _1(x))\subseteq S_1(x)$, 
then by continuity of $\alpha _1$ (and the fact that $x\in B_1$) there is 
a definable open neighbourhood $B$ of $x$ such that
$\alpha _1(B)\cap (z',\alpha _1(x))$ is infinite. But then, since
$\alpha _1(Ox\cap B)\cap (z',\alpha _1(x))$ is 
infinite (because, otherwise we would have $x\in A_1$), 
$T_1(x)$ is infinite as well. 
Since $O$ has strong definable choice, we have a strong definable 
choice $l'_1$ for the definable family $\{T_1(x):x\in B_1\}$ and from 
this we get $l_1$ for the definable family
$\{V_0(x):x\in B_1\}$ by setting $l_1(x):=\alpha _1(l'_1(x)x)$. Note
that if $V_0(x)=V_0(y)$ then $V_0(l'_1(x)x)=V_0(l'_1(y)y)$.
Let $U_1:=X_1\cup Y_1\cup M_1$ then $U_1$ is large in $U_0$ and
the claim is proved for $k=1$.

Suppose that the claim is true for $k$. We will show that it is 
true for $k+1$.
For this consider the definable family $\{V_k(x):x\in U_k\}$
of non empty definable subsets of $U$, where
$V_k(x):=\{u\in xV:\pi _k(u)=(l_1(x),\dots ,l_k(x))\}$
(note that we have $xV=yV$ iff $V_k(x)=V_k(y)$), and substitute in the
proof for the case $k=1$, $0$ by $k$ and $1$ by $k+1$.
\qed

\end{section}

\begin{section}{Definable extensions}
\label{section definable extensions}

\begin{subsection}{Definable $G$-modules}
\label{subsection definable $G$-modules}

\begin{defn}\label{defn G-module}
{\em Let $G$ be a definable group. A {\it definable} $G$-{\it module}
is a pair $(A, \gamma )$ where $A$ is a definable abelian group
and $\gamma :G\into Aut_{\N}$$(A)$ is a homomorphism form $G$ into
the group of all definable automorphisms of $A$, such that the map
$\gamma :G\times A\into A,$ $\gamma (x,a):=\gamma (x)(a)$ is
definable.

We say that $A$ is {\it trivial} if 
$\forall x\in G \forall a\in B,\, \gamma (x)(a)=a$, 
$A$ is {\it faithful} if $\gamma :G\into Aut_{\N}$$(A)$ is injective. 
A definable {\it $G$-submodule} of $A$ is a definable normal subgroup
$B$ of $A$ such that $\forall x\in G,\, \gamma (x)(B)\subseteq B$ 
(i.e., $B$ is invariant under $\gamma $). We then
have natural induced definable $G$-modules $(B,\gamma _{|B})$
and $(A/B,\gamma _{A/B})$. 
We say that $A$ is {\it irreducible} if it has no proper
definable $G$-submodules. A special definable $G$-submodule of $A$ is
$A^G:=\{ a\in A:\forall x\in G,\, \gamma (x)(a)=a\}$.
}
\end{defn}

The next lemma follows from theorem \ref{thm definable quotients} but
we include here a direct prove based on DCC.

\begin{lem}\label{lem quotient by AG and ker are definable}
Let $(A, \gamma )$ be a definable $G$-module. Then $A/A^G$ is a
definable group,
$Ker \gamma $ is a normal definable subgroup of $G$,
$\bar{G}:=G/Ker \gamma $ is definable and we
have a natural induced faithful definable $\bar{G}$-module
$(\bar{\gamma }, A)$. 

Also, if $U$ is a definable group and $A$ is a normal
subgroup of $U$ then $C_U(A)$ is a normal definable subgroup of $U$ 
and $U/C_U(A)$ is definable. In particular, $U/Z(U)$ is definable.
\end{lem}

\pf
For each $g\in G$ we have a definable endomorphism
$\alpha (g):A\into A$ given by
$\forall a\in A,\, \alpha (g)(a):=\gamma (g)(a)-a$
and
$A^G=\bigcap _{g\in G}ker\alpha (g)$. And so by DCC on definable
subgroups (see \cite{p}) there are
$g_1,\dots, g_n\in G$ such that $A^G=\bigcap _{i=1}^n
ker\alpha (g_i)$. But then, the definable map $A\into 
\alpha (g_1)(A)\times \cdots \times \alpha (g_n)(A)$, 
$a\into (\alpha (g_1)(a),\dots ,\alpha (g_n)(a))$ shows that $A/A^G$
is definable.

Let $a\in A$ and consider the definable map 
$\beta (a):G\into A$, $g\into \gamma (g)(a)-a$ then $\{g\in G:\beta
(a)(g)=0\}$ is a definable subgroup of $G$ and $Ker \gamma =\bigcap
_{a\in A}\{g\in G:\beta (a)(g)=0\}$ and by DCC on definable subgroups
there are $a_1,\dots ,a_n\in A$ such that $Ker \gamma =\bigcap_{i=1}^n
\{g\in G:\beta (a_i)(g)=0\}$. The definable map $G\into 
\beta (a_1)(G)\times \cdots \times \beta (a_n)(G)$, 
$g\into (\beta (a_1)(g),\dots ,\beta (a_n)(g))$ shows that 
$G/Ker\gamma $ is definable. 

If $U$ is a definable group and $A$ is a normal
subgroup then $C_U(A)=\bigcap _{a\in A}C_U(a)$ and by DCC on definable 
subgroups there are $a_1,\dots ,a_n\in A$ such that 
$C_U(A)=\bigcap _{i=1}^nC_U(a_i)$ and so $C_U(A)$ is definable (and
normal) and if for each $a\in A$ we define $ad(a):U\into U$ by
$\forall u\in U,\, ad(a)(u):=aua^{-1}u^{-1}$ then  
the definable map $U\into ad(a_1)(U)\times \cdots \times
ad(a_n)(U)$, $u\into (ad(a_1)(u),\cdots ,ad(a_n)(u))$
shows that $U/C_U(A)$ is definable.
\qed

\end{subsection}

\begin{subsection}{Group cohomology}
\label{subsection group cohomology}

For the rest of this subsection we assume that
$(A,\gamma )$ is a definable $G$-module.

\begin{defn}\label{defn n-cochain}
{\em For each $n\in \Bbb N$ let $C^n_{\N}$$(G,A,\gamma )$ denote the 
abelian group of all definable
functions from $G^n$ into $A$ with point wise addition.
An element of $C^n_{\N}$$(G,A,\gamma )$ is called a {\it definable} 
$n$-{\it cochain (over $\N$}).}
\end{defn}

\begin{defn}\label{defn co-boundary map}
{\em The {\it co-boundary map} $\delta :C^n_{\N}$$(G,A,\gamma )\into
C^{n+1}_{\N}$$(G,A,\gamma ),$
is defined by 
\[  \delta (c)(g_1,\dots ,g_{n+1}):=\gamma (g_1)(c(g_2,\dots ,g_{n+1})) +\] 
\[ \sum_{i=1}^n(-1)^ic(g_1,\dots ,g_ig_{i+1},\dots ,g_{n+1}) 
 + (-1)^{n+1}c(g_1,\dots ,g_n) .\]
}
\end{defn}

It is clear that $\delta(c)$ is also definable.

\begin{lem}\label{lem co-boundary map}
$\delta \delta =0.$
\end{lem}

\pf
This is a simple calculation.
\qed

\begin{defn}\label{defn n-cohomology group}
{\em We therefore have a {\it complex} $C^*_{\N}$$(G,A,\gamma ).$ 
$ B^n_{\N}$$(G,A,\gamma )$ denotes the image of
$\delta :C^{n-1}_{\N}$$(G,A,\gamma
)\into C^n_{\N}$$(G,A,\gamma ),$
$ Z^n_{\N}$$(G,A,\gamma )$ denotes the kernel of
$\delta :C^n_{\N}$$(G,A,\gamma )\into 
C^{n+1}_{\N}$$(G,A,\gamma ),$
and 
$ H^n_{\N}$$(G,A,\gamma )$ denotes
$Z^n_{\N}$$(G,A,\gamma )/B^n_{\N}$$(G,A,\gamma ).$
$H^n_{\N}$$(G,A,\gamma )$ is the $n$-{\it cohomology group over $\N$}, 
the elements of
$B^n_{\N}$$(G,A,\gamma )$ are the {\it definable} $n$-{\it coboundaries} and  
the elements of
$Z^n_{\N}$$(G,A,\gamma )$ are the {\it definable} $n$-{\it cocycles}.
}
\end{defn}

\begin{nrmk}\label{nrmk A1 times A2}
{\em
Let $(A,\gamma )$ be a definable $G$-module. 
Suppose that $A:=A_1\times A_2$
and that $A_1$ and $A_2$ are invariant under the action
of $G$ on $A$.
Then $H^n_{\N}$$(G,A,\gamma )$ is isomorphic with 
$H^n_{\N}$$(G,A_1,\gamma _{|A_1})
\times H^n_{\N}$$(G,A_2,\gamma _{|A_2})$.
}
\end{nrmk}

\end{subsection}

\begin{subsection}{Definable extensions}
\label{subsection definable extensions}

\begin{defn}\label{defn extension and section}
{\em Let $U$ be a definable group. $(U,i,j)$ is an {\it definable extension} of
$G$ by $A$ if we have an exact sequence
\[ 1\rightarrow A\stackrel{i}{\rightarrow}U\stackrel{j}{\rightarrow}
G\rightarrow 1 \]
in the category of definable groups with definable homomorphisms.
A {\it definable section} is a definable map $s:G\into U$
such that $\forall g\in G, \, j(s(g))=g$.
}
\end{defn}

\medskip
\textbf{Note:} Below we will some times assume that $A\unlhd U$, and write 
$(U,j)$ for $(U,i,j)$.

\begin{thm}\label{thm existence of definable sections}
Let 
$1\rightarrow A\rightarrow U\stackrel{j}\rightarrow G\rightarrow 1$
be a definable extension. Then there is a definable section
$s:G\into U$.
\end{thm}

\pf
Suppose that 
$U\subseteq N^m$. For each $l\in \{1,\dots , m\}$ let 
$\pi _l:N^m\into N^l$ be the projection into the first $l$ 
coordinates and let $\pi ^l:N^m\into N$ be the projection onto the
$l$-th coordinate. The proof of the theorem follows from the proof
of theorem \ref{thm definable quotients} after making the 
following substitutions: $U_0:=G$,
for each $x\in U_0$, $V_0(x):=j^{-1}(x)$ and the definable 
neighbourhoods in $U$ that appear in the proof of theorem 
\ref{thm definable quotients} are substituted by definable 
neighbourhoods in $G$.
\qed

\begin{defn}\label{defn equivalent extensions}
{\em Two definable extensions
$1\rightarrow A\stackrel{i}{\rightarrow}U\stackrel{j}{\rightarrow}
G\rightarrow 1$ and 
$1\rightarrow A\stackrel{i'}{\rightarrow}U'\stackrel{j'}{\rightarrow}
G\rightarrow 1$ are
{\it definably equivalent}
if there is a definable
homomorphism $\varphi: U\into U'$ such that 

\[ 
\begin{array}{clcr}
                                         U \\ 
               \stackrel{i}{\nearrow}\;\;\;\; \stackrel{j}{\searrow}  \\
1\rightarrow A\;\;\;\;\;\stackrel{\varphi}{\downarrow}\;\;\;\;\;G\rightarrow 1\\
             \stackrel{i'}{\searrow}\;\;\;\;\stackrel{j'}{\nearrow} \\
                                        U' \\
\end{array}
\]

is a commutative diagram.
}
\end{defn}
\end{subsection}

\begin{subsection}{Definable $G$-kernels}
\label{subsection definable $G$-kernels}

\textbf{Notation:} Let $A$ be a definable group. 
$Aut_{\N}$$(A)$ denotes the group of all definable automorphisms of $A$, 
$Inn(A)$ the group of all inner automorphisms of $A$ and 
$Out_{\N}$$(A):=Aut_{\N}$$(A)/Inn(A)$. Let $\iota
:Aut_{\N}(A)\into Out_{\N}$$(A)$ denote the natural homomorphism.
If $A\unlhd U$ and $u\in U$ then we denote by $<u>$ the 
automorphism of $A$ given by 
$<u>(a):=uau^{-1}$ for all $a\in A$.

\begin{defn}\label{defn G-kernel}
{\em
Let $G$ be a definable group. A {\it definable $G$-kernel} 
$(A,\theta )$ is a definable group $A$ with a homomorphism
$\theta :G\into Out_{\N}$$(A)$ and a homomorphism 
$\alpha :G\into Aut_{\N}$$(A)$ such that  $\theta =\iota \circ
\alpha $ and the map $\alpha :G\times A\into A$, $\alpha (g,a):=\alpha
(g)(a)$ is definable.
Note that $\theta $ induces a definable action $\theta _0 :G\times
Z(A)\into Z(A)$ making the center $Z(A)$ of $A$ a {\it definable
$G$-module}. 
We say that $\alpha $ as above is a {\it definable representative} of
the definable $G$-kernel $(A,\theta )$ and we write $\alpha \in \theta
$. 
}
\end{defn}

If $\alpha ,\beta \in \theta $ then by theorem \ref{thm definable choice}
there is a definable function $k:G\into A$ such that 
$\forall x\in G,\, \beta (x)=<k(x)>\alpha(x)$. 
Note also that, by theorem \ref{thm definable choice}
there is a definable function $h_{\alpha }:G\times G\into A$
such that we have $\forall x,y\in G,\, h_{\alpha }(x,1)=h_{\alpha
}(1,y)=1$, and
\begin{eqnarray}\label{eq h and alpha}
\forall x,y\in G,\, \alpha (x)\alpha (y)=<h_{\alpha }(x,y)>\alpha (xy)
\end{eqnarray}
and 
$\forall x,y\in G,\, \beta (x)\beta (y)=<h_{\beta }(x,y)>\beta (xy)$ where
$h_{\beta }:G\times G\into A$ is the definable function given by 
\[ \forall x,y\in G,\, 
h_{\beta }(x,y):=k(x)\alpha (x)(k(y))h_{\alpha }(x,y)k(xy)^{-1}.\]

Note that if $(A,\theta )\in EK_{\N}$$(G,B)$ and let 
$(U,\pi )$ is
a definable extension of $A$ by $G$ and  $s:G\into U$ is a definable
section. Then 
\[ \forall x,y\in G,\, h_{\alpha _{U,s}}(x,y):=s(x)s(y)s(xy)^{-1}.\]

\begin{defn}\label{defn group of kernels}
{\em
Let $G$ be a definable group and $B$ an abelian definable group.
Two definable $G$-kernels $(A_i,\theta _i)$ with $i=1,2$ with centre
$B$ are {\it definably equivalent} if there is a definable isomorphism
$\sigma :A_1\into A_2$ such that for all $b\in B$ and for all $x\in G$,
$\sigma (b)=b$ and
$\sigma \theta _1(x)\sigma ^{-1}= \theta _2(x)$.
This relation is an equivalence relation and the set of all the 
classes is denoted by $K_{\N}$$(G,B)$. 
}
\end{defn}

\begin{nrmk}\label{nrmk canonical action}
{\em Let $(U,\pi)$ be a definable extension of $G$ by $A$. Then 
there is a canonical homomorphism $\theta _{U}:G\into Out_{\N}$$(A)$
such that 
$(A,\theta _U)$ is a definable $G$-kernel: take, for each $x\in G$,
$\theta _U(x):=\{ <u>:u\in \pi ^{-1}(x)\}$ with definable
representative given by $\alpha _{U,s}:G\into Aut_{\N}$$(A)$, 
$\alpha _{U,s}(g)(a):=<s(g)>(a)$ where $s:G\into U$ is a definable section.
}
\end{nrmk}

\begin{defn}\label{defn EK}
{\em
A definable $G$-kernel $(A,\theta )$ is {\it definably extendible}
if there is a definable extension $(U,\pi )$ of $G$ by $A$ such 
that $(A,\theta _U)$ is definably equivalent to $(A,\theta )$.
We say in
this case that $(U,\pi)$ is {\it compatible} with the $G$-kernel.
We denote by $Ext_{\N}$$(G,A,\theta )$ the set of all equivalence classes of
definable extensions of $G$ by $A$ compatible with the $G$-kernel 
$(A, \theta)$. 
Let $EK_{\N}$$(G,B)$ be the subset of $K_{\N}$$(G,B)$ of all 
classes $(A,\theta )$ such that $Ext_{\N}$$(G,A,\theta )$ is nonempty.
Note that $EK_{\N}$$(G,B)$ is a well defined subset of 
$K_{\N}$$(G,B)$.}
\end{defn}

\end{subsection}

\begin{subsection}{Existence of definable extensions}
\label{subsection existence of definable extensions}

With the set up we have established, the proof of the following results
is now as in the classical case, 
for details see the relevant lemmas in \cite{em1} and
\cite{em2}.

\begin{fact}\label{fact main fact in em}
{\em
There is a canonical map from
$K_{\N}$$(G,B)$ into $H^3_{\N}$$(G,B,\theta _0)$, sending
$(A,\theta )$ into $c_{(A,\theta )}$
and $(A,\theta )\in EK_{\N}$$(G,B)$ if and only if 
$c_{(A,\theta )}=1$.

Let $(A, \theta )\in K_{\N}$$(G,B)$ and let $\alpha \in \theta $ and 
let $h_{\alpha }$ be the corresponding definable function as in 
equation \ref{eq h and alpha}. 
For $x,y,z\in G$, using associativity, the product
$\alpha (x)\alpha (y)\alpha (z)$ may be calculated in two different 
ways. The identity of this two results gives  
for all $x,y,z\in G$, the following identity
$<h_{\alpha }(x,y)h_{\alpha }(xy,z)>=$
$<\alpha (x)(h_{\alpha }(y,z))h_{\alpha }(x,yz)>$.

But only the elements of the center $B$ of $A$ determine the identity
inner automorphism. Hence there exists a definable $3$-cochain
$c_{\alpha }
\in C^3_{\N}$$(G,B,\theta _0)$ such that 
\begin{eqnarray}\label{eq h and c}
\forall x,y,z\in G,\, \alpha (x)(h_{\alpha }(y,z))h_{\alpha
}(x,yz)=c_{\alpha }(x,y,z)h_{\alpha }(x,y)h_{\alpha }(xy,z).
\end{eqnarray}

Now some calculations show that 
$c_{\alpha }\in Z^3_{\N}$$(G,B,\theta _0)$ and if $\beta \in \theta $ then 
$h_{\beta }(x,y)=g_{\alpha ,\beta }(x,y)h_{\alpha }(x,y)$
where $g_{\alpha ,\beta }\in C^2_{\N}$$(G,B,\theta _0)$) and
$c_{\alpha }$
is changed to a cohomologous cocycle $c_{\beta }$ and 
by suitably changing the choice of $\alpha \in \theta $, $c_{\alpha }$ 
may be changed to any cohomologous cocycle.

Suppose now that $(A,\theta )\in EK_{\N}$$(G,B)$ and let 
$(U,\pi )$ be
a definable extension of $A$ by $G$ and let $s:G\into U$ be a definable
section. Then a simple calculation
shows that 
\begin{eqnarray}\label{eq h for U s}
\alpha _{U,s}(x)(h_{\alpha _{U,s}}(y,z))h_{\alpha _{U,s}}(x,yz)=
h_{\alpha _{U,s}}(x,y)h_{\alpha _{U,s}}(xy,z)
\end{eqnarray}
and therefore $c_{\alpha _{U,s}}(x,y,z)=1$. 

Conversely, suppose that $(A,\theta )\in K_{\N}$$(G,B)$ is such that 
$c_{(A,\theta )}=1$ in $H^3_{\N}$$(G,B,\theta _0)$. 
Select $\alpha \in \theta $ such that $c_{\alpha }(x,y,z)=1$ 
for all $x,y,z\in G$. The proof of the 
result below shows that we can find $(U,\pi )\in 
Ext_{\N}$$(G,A,\theta )$.
}
\end{fact}

\begin{fact}\label{fact main fact}
{\em
Let $(A,\theta )\in EK_{\N}$$(G,B)$ and $(U,\pi )\in 
Ext_{\N}$$(G,A,\theta )$. Then
there is a canonical bijection from $Ext_{\N}$$(G,A,\theta )$
into $H^2_{\N}$$(G,B,\theta _0)$ sending $(U,\pi )$ into the identity of 
$H^2_{\N}$$(G,B,\theta _0)$.

Let $s:G\into U$ be a definable section. We construct 
$(V_s,i_s,j_s)\in Ext_{\N}$$(G,A,\theta )$ associated
with $s$ as follows: $V_s$ as domain $A\times G$ and multiplication given by
\begin{eqnarray}\label{eq group of h} 
\forall a,b\in A\forall x,y\in G,\,(a,x)(b,y)=
(a[\alpha _{U,s}(x)(b)]h_{\alpha _{U,s}}(x,y),xy).
\end{eqnarray}
From equation (\ref{eq h for U s}) and equation (\ref{eq h and c})
$V_s$ is a definable group, $(1,1)$ is the
identity, and the inverse of $(a,x)$ is 
$(\alpha _{U,s}(x)^{-1}[h_{\alpha _{U,s}}(x,x^{-1})a]^{-1},x^{-1})$, 
$i_s:A\into V_s$ given by
$\forall a\in A,\, i(a):=(a,1)$ and $j_s :V_s\into G$ given by
$\forall a\in A\forall x\in G,\,j(a,x):=x$. The map $t:G\into V_s$ given
by $\forall x\in G,\, t(x):=(1,x)$
is a definable section, we see that for all $x\in G,\, <t(x)>=
\alpha _{U,s}(x)$ and therefore 
$(V_s,i_s,j_s)\in Ext_{\N}$$(G,A,\theta )$.
Also, the map $U\into V_s$, $u:=as(x)\into (a,x)$ is a definable isomorphism.

Moreover, if $s':G\into U$ is another definable section and 
$(V_{s'},i_{s'},j_{s'})\in Ext_{\N}$$(G,A,\theta )$ 
the corresponding
definable extension, then there is a definable function $k_{s,s'}:G\into A$
given by $\forall x\in G,\, s'(x):=k_{s,s'}(x)s(x)$ such that 
\begin{eqnarray}\label{eq h another section} 
\forall x,y\in G,\, h_{\alpha _{U,s'}}(x,y):=
k_{s,s'}(x)\alpha _{U,s}(x)(k_{s,s'}(y))h_{\alpha _{U,s}}(x,y)
k_{s,s'}(xy)^{-1}
\end{eqnarray}
and the map $V_s\into V_{s'}$, $(a,x)\into (ak_{s,s'}(x)^{-1},x)$ is a 
definable isomorphism.

Also, $V_s$ (and therefore $U$) is definably isomorphic with 
$A\rtimes _{\alpha _{U,s}}G$ iff
there is a definable function $g:G\into A$ such that 
\begin{eqnarray}\label{eq semi direct} 
\forall x,y\in G,\, h_{\alpha _{U,s}}(x,y)=\alpha _{U,s}(x)
(g(y))g(x)g(xy)^{-1},
\end{eqnarray}
since if $g:G\into A$ satisfying equation (\ref{eq semi direct})
, then the function $G\into V_s$, $x\into (g(x)^{-1},x)$ is a
homomorphism.

Finally, if $(U',\pi' )\in Ext_{\N}$$(G,A,\theta )$ 
and $s':G\into U'$ is a definable section and
$h_{\alpha _{U',s'}}:G\times G\into A$ is the corresponding definable 
function then
there is $c\in Z^2_{\N}$$(G,B,\theta _0)$ such that
\begin{eqnarray}\label{eq h another extension} 
\forall x,y\in G,\, h_{\alpha _{U',s'}}(x,y)=c(x,y)h_{\alpha _{U,s}}(x,y),
\end{eqnarray} 
$c$ in $H^2_{\N}$$(G,B,\theta _0)$ does not depend on the equivalence class
of $(U',\pi' )$ or the on the choice of the definable
section. Moreover, $c$ is zero in $H^2_{\N}$$(G,B,\theta _0)$ iff
$(U,\pi )$ and $(U',\pi' )$ are definably equivalent. 
}
\end{fact}

\begin{fact}\label{fact extension 2-co-cycle}
{\em
Let $(A,\gamma )$ be a definable $G$-module. Then there is a bijection
from $Ext_{\N}$$(G,A,\gamma )$ onto $H^2_{\N}$$(G,A,\gamma )$ 
sending the class of $A\rtimes _{\gamma }G$ into the identity of 
$H^2_{\N}$$(G,A,\gamma )$.

Let $(U,\pi)\in Ext_{\N}$$(G,A,\gamma )$ and let $s:G\into U$ be a 
definable 
section. Then there is a
canonical definable $2$-cocycle $c\in Z^2_{\N}$$(G,A,\gamma )$ 
associated with this 
definable section given by: $\forall g,h\in G,\, 
c(g,h):=s(g)s(h)s(gh)^{-1}$, and therefore, (in $A$) we have

\begin{eqnarray}\label{eq c(g,h) in Z2}
\forall g,h,k \in G, \,c(h,k)^g-c(gh,k)+c(g,hk)-c(g,h)=0,
\end{eqnarray}
and, there is $(V,i,j)\in Ext_{\N}$$(G,A,\gamma )$
associated with the definable $2$-cocycle $c$ given by:
$V:=A\times G$ 
and with  multiplication given by

\begin{eqnarray}\label{eq group for c(g,h)}
\forall a,b\in A, \,\forall g,h \in G, \,(a,g)(b,h):=(a+b^g+c(g,h),gh),
\end{eqnarray} 
from equation 
(\ref{eq c(g,h) in Z2}) $V$ is a group, with identity $(-c(1,1),1),$
$i:A\into V$ given by $i(a):=(a-c(1,1),1)$ and $j:V\into G$ by 
$j(a,g):=g.$ And the map $U\into V$, $u:=as(g)\into (a,g)$ is a definable 
isomorphism.

If $s':G\into U$ is
another definable section, and $c'\in Z^2_{\N}$$(G,A,\gamma )$ is the 
corresponding 
definable $2$-cocycle and  $(V',i',j')\in Ext_{\N}$$(G,A,\gamma )$
is the corresponding definable extension, then 
there is a definable function $b:G\into A$ given by, $\forall g\in G,\,
s'(g):=b(g)s(g)$ such that 

\begin{eqnarray}\label{eq c(g,h) and c'(g,h)}
\forall g,h \in G, \, c'(g,h)-c(g,h)=b(h)^g-b(gh)+b(g),  
\end{eqnarray}
i.e., $c$ and $c'$ determine the same element in
$H^2_{\N}$$(G,A,\gamma )$ and the map 
$V\into V'$, $(a,g)\into (a-b(g),g)$ is a definable isomorphism. 

Also, $V$ (and therefore $U$) is definably isomorphic with 
$A\rtimes _{\gamma }G$ iff 
there is a definable function $a:G\into A$ such that  

\begin{eqnarray}\label{eq c(g,h) split}
\forall g,h \in G, \, c(g,h)=a(h)^g-a(gh)+a(g),  
\end{eqnarray}
(since if $a:A\into G$ exists and satisfies equation (\ref{eq c(g,h) split}) 
, the function $G\into V$, $g\into (-a(g),g)$ is 
a homomorphism). 
}
\end{fact}

\begin{fact}\label{fact groups on ext}
{\em
The set $K_{\N}$$(G,B)$ can be made into an abelian
group in the following
way: Let $(A_i,\theta _i)\in K_{\N}$$(G,B)$ with $i=1,2$. Then their
{\it product} $(A_1,\theta _1)\otimes (A_2,\theta _2)
\in K_{\N}$$(G,B)$ is the element, which is the class of 
$(A,\theta )$ where
$A:=A_1\times A_2/C$ where $C:=\{(b,b^{-1}):b\in B\}$, and
$\theta $ is represented by the map 
$\alpha:=(\alpha _1,
\alpha _2):G\times A_1\times A_2\into A_1\times A_2$, given by
$\alpha (g)(a_1,a_2):=(\alpha _1(g)(a_1),\alpha _2(g)(a_2))$, where
$\alpha _i$ represents $(A_i,\theta _i)$.
The identity of $K_{\N}$$(G,B)$ is the definable $G$-kernel 
$(B,\theta _0)$. 
And the inverse of $(A,\theta)\in K_{\N}$$(G,B)$ is $(A^{*},\theta)$
where $A^{*}$ is the anti-isomorphic with $A$ with domain $A$.

$EK_{\N}$$(G,B)$ is a subgroup of $K_{\N}$$(G,B)$: 
Let $(A_i,\theta _i)\in 
EK_{\N}$$(G,B)$ with definable extensions $(U_i,\pi _i)$$
\in Ext_{\N}$$(G,A_i,\theta _i)$ where $i=1,2$. Then, if $(A,\theta):=
(A_1,\theta _1)\otimes (A_2,\theta _2)$ then $(U,\pi):=
(U_1,\pi _1)\otimes (U_2,\pi _2)$ is an element of
$Ext_{\N}$$(G,A,\theta )$, where 
$U:=D/E$, $D:=\{(u_1,u_2):u_1\in U_1,u_2\in U_2,
\pi _1(u_1)=\pi _2(u_2)\}$, $E:=\{(b,b^{-1}):b\in B\}$ and $\pi $
is induced by any of $\pi _i$.
If
$(A,\theta )\in EK_{\N}$$(G,B)$ with $(U,\pi)\in Ext_{\N}$$(G,A,\theta)$ then
$(U^*,\pi ^*)\in Ext_{\N}$$(G,A^*,\theta)$ where $U^*$ is the group anti-
isomorphic with $U$ with domain $U$ and for all $u\in U^*$, 
$\pi ^*(u):=\pi (u^{-1})$.
The group $K_{\N}$$(G,B)/EK_{\N}$$(G,B)$ is called
the group of {\it similarity} classes. Note that 
$Ext_{\N}$$(G,B, \theta _0)$ can be made into a group with product 
$\otimes $ defined above.

The map given in fact \ref{fact extension 2-co-cycle} is a
isomorphisms between $H^2_{\N}$$(G,B, \theta _0)$ and $Ext_{\N}$
$(G,B,\theta _0)$. Moreover, the map from
$H^2_{\N}$$(G,B,\theta _0)$
into $Ext_{\N}$$(G,A,\theta )$ for a fixed $(U,\pi )$ in
$Ext_{\N}$$(G,B,\theta )$
of fact \ref{fact main fact in em} is 
the composition of the isomorphism from
$H^2_{\N}$$(G,B,\theta _0)$ into $Ext_{\N}$$(G,B,\theta _0)$ 
and the map from $Ext_{\N}$$(G,B,\theta _0)$
into $Ext_{\N}$$(G,A,\theta )$ which sends
$(V,j)$ into $(U,\pi )\otimes (V,j)$.
Finally note that the map from fact \ref{fact main fact} is a 
homomorphism with kernel $EK_{\N}$$(G,B)$.
}
\end{fact}

\end{subsection}
\end{section}

\begin{section}{Definably compact definable groups}

In this section we prove that a definably compact definable group
is abelian-by-finite. This will follow after we show that a definably
connected definably compact definable $G$-module where $G$ is infinite
and definably connected is trivial.
Before we proceed, we need the following easy lemma.

\begin{lem}\label{lem nice curve}
Let $U$ be an infinite  definable group and let $V$ be a definable
subgroup such that $dimV<dimU$. Then there is a definable continuous 
embedding $\sigma :(a,b)\into U$
such that $\lim _{t\into a^{+}}\sigma (t)=0$ and $\sigma
(a,b)\subseteq U\setminus V$.
\end{lem}

\pf
Let $(O,\phi )$ be a definable chart of 1 (the identity  of
$U$). Then $\phi (O)$ is a definable open subset of $N^n$ where
$n=dimU$. Let $e=(e_1,\dots ,e_n)=\phi (1)$, $B=I_1\times
\cdots \times I_n\subseteq \phi (O)$ an open box containing $e$ and for each
$i=1,\dots ,n$ let $\bar{I_i}:=\{e_1\}\times \cdots \times
\{e_{i-1}\}\times I_i\times \{e_{i+1}\}\times \cdots \times \{e_n\}$
and $J_i:=\phi ^{-1}(\bar{I_i})$. Let $D:=\phi (V\cap O)\cap B$. If
there is $i\in \{1,\dots ,n\}$ and an open subinterval $I$ of
$\bar{I_i}$ with one endpoint $e$ and such that $I\cap D=\emptyset $
then we are done. So suppose otherwise. Then after substituting each
$I_i$ with a smaller interval if necessary, 
we have that each $J_i\subseteq V$. But
this clearly implies that $dimV=n$.  
\qed

\begin{thm}\label{thm compact G modules}
Let $(A,\gamma )$ be a definably compact, definably connected 
definable $G$-module, where $G$ is an infinite definably
connected definable group. Then $(A,\gamma )$ is trivial.
\end{thm}

\pf
Without loss of generality, we can assume that $\N$ is $\aleph
_1$-saturated, in particular $|N|>\aleph _0$.
Suppose that $A^G\neq A$, and let $B$ be an infinite minimal definable
subgroup of $A/A^G$. Let $C$ be a definable subgroup of $A$ such that
$C/A^G=B$ and let $\bar{C}$ be the smallest definable $G$-submodule of
$A$ containing $C$. Then  $C\in Ext_{\N}$$(B,A^G)$ 
i.e., $C$ is a definable extension of $B$ by $A^G$ and
there is a definable section $s:B\into C$.
Let $c(b,b'):=s(b)+s(b')-s(b+b')$
be the corresponding definable $2$-cocycle, then we have a
definable family $\Gamma :G\times B\into \bar{C}$ of definable homomorphisms
from $B$ into $\bar{C}$ given by, $\forall
g\in G\forall b\in B,\, \Gamma (g,b)=\gamma (g)(s(b))-s(b)$ and such
that $\forall g\in G\forall b\in B,\, \Gamma (1,b)=0=\Gamma (g,0)$. 
To see this, subtract to the equation above for the $2$-cocycle the
equation obtained from it after applying $\gamma (g)$.
Since for
each $c\in C$ there are unique $a\in A^G$ and $b\in B$ such that
$c=a+s(b)$ and for all $g\in G$, $\gamma (g)(c)=a+\gamma (g)(s(b))$
we must have $ker_G\Gamma \neq G$ where $ker_G\Gamma :=\{g\in
G:\forall b\in B,\, \Gamma (g,b)=0\}$.

Since $B$ has no infinite proper definable subgroups,
for each $g\in G$, $\Gamma (g)(B)$ is either $0$ or infinite (with the
same dimension as $B$) and with no infinite proper definable additive
subgroups and so 
by dimension consideration, there is a minimal $n\geq 1$ such that 
for each $i\in \{1,\dots ,n\}$ there is $g_i\in G$
such that $\Gamma (g_i)(B)\neq 0$ and
$\Gamma (G)(B)\subseteq \Gamma (g_1)(B)+\cdots +\Gamma (g_n)(B)$. Now since
$F:=\bigcap _{i=1}^n\Gamma (g_i)(B)$ is finite, $D:=\Gamma (G)(B)/F$
is definable and we have a natural induced definable family
$\Lambda :G\times B\into D$ of 
definable homomorphisms from $B$ into $D$. 
Its easy to see that $ker_G\Lambda \neq G$.
Now for each $i\in \{1,\dots ,n\}$ let $D_i:=\Gamma (g_i)(B)/F$. Then 
$D=\bigoplus _{i=1}^nD_i$ and we have natural induced
definable families $\Lambda _i:G\times B\into D_i$
of definable homomorphisms from $B$ into $D_i$, and there is 
$i_0\in \{1,\dots ,n\}$ such that $ker_G\Lambda _{i_0}\neq G$.

Since for each $g\in G\setminus ker_G\Lambda _{i_0}$, 
$ker\Lambda _{i_0}(g)$ is finite, and $Tor(B)$ is countable (see
corollary 5.8 \cite{s}) therefore 
$\cup \{ker\Lambda _{i_0}(g):g\in G\setminus ker_G\Lambda _{i_0}\}$
is finite and so there is a finite additive subgroup
$E$ of $B$ such that for all $g\in G\setminus ker_G\Lambda _{i_0}$, 
$ker\Lambda _{i_0}(g)\subseteq E$. Let $B':=B/E$. Its easy to see that 
$\Lambda _{i_0}$ induces a natural definable family
$\Phi :G\times B'\into B'$ of definable endomorphisms of $B'$
such that $ker_G\Phi \neq G$ and 
for each $g\in G\setminus ker_G\Phi $, $\Phi (g)$ is a definable
automorphism of $B'$.

Since $ker_G\Phi \neq G$ and $G$ is definably connected,
we have $dim(ker_G\Phi )<dimG$ 
and by lemma \ref{lem nice curve}  
there is a definable continuous embedding $\sigma :(a,b)\into G$
such that $\lim _{t\into a^{+}}\sigma (t)=1$ and $\sigma
(a,b)\subseteq G\setminus ker_G\Phi $. Let $x_0\in B'\setminus
\{0\}$. Then for every $t\in (a,b)$ there exists a unique $x\in B'$
such that $\Phi (\sigma (t),x)=x_0$. This gives us a definable
function $\tau :(a,b)\into \tau (a,b)\subseteq B'$. 
Since $B'$ is definably compact, there is
an element $c\in B'$ such that $\lim _{t\into a^{+}}\tau (t)=c$. But then, by 
continuity of $\Phi $ we have $0=\Phi (1,c)=x_0$, 
and so we get a contradiction. 
\qed

The next corollary was also proved in \cite{pst2} but assuming that
$\N$ has definable Skolem functions and using the theory of $\bigvee
$-definable groups.

\begin{cor}\label{cor compact U}
Let $U$ be a definably compact, definably connected
definable group. Then $U$ is either abelian or
$U/Z(U)$ is a definable semi-simple group. 
In particular, if $U$ is solvable then it is abelian. 
\end{cor}

\pf
By lemma \ref{lem quotient by AG and ker are definable}, $U/Z(U)$ is definable.
Suppose that $U/Z(U)$ is infinite and not semi-simple. Then there is a
normal definably connected definable subgroup $X$ of $U$ such that $Z(U)\leq X$
and $X/Z(U)$ is an abelian infinite normal definable subgroup of $U/Z(U)$. 
Now $X$ is a definable $U$-module by conjugation and by
theorem \ref{thm compact G modules}, $X=X^U\leq Z(U)$ contradiction. 
\qed

\end{section}

\begin{section}{Definable solvable groups}
\label{section definable solvable groups}

\begin{subsection}{Preliminary lemmas}
\label{subsection preliminary lemmas}

The results below are stated for definable $G$-modules, but 
each one of them has a corresponding analogue for definable $G$-kernels.
These are obtained after making the obvious substitutions.
Since after these substitutions, the proofs are exactly the same, we
omit them. We will be using through this subsection the results of
subsection \ref{subsection existence of definable extensions}. We will
also often use the following fact: 

\begin{fact}\label{fact Al split}
{\em
Let $A:=A_1\times \cdots \times A_k$ and suppose that $(A,\gamma )$
is a definable $G$-module, and let 
$(U,j)\in Ext_{\N}$$(G,A,\gamma )$ with the corresponding
canonical definable $2$-cocycle
$c\in Z^2_{\N}$$(G,A,\gamma )$. Suppose also that each $A_i$ is
invariant under $G$, then 
$c:=(c_1, \dots ,c_k)$ where for each $i\in \{1,\dots ,k\}$, 
$c_i\in Z^2_{\N}$$(G,A_i,\gamma _{|A_i})$. Let $l\in \{1,\dots ,k\}$. 
If for each
$i\in \{1,\dots ,k\}\setminus \{l\}$, $c_i\in B^2_{\N}$
$(G,A_i,\gamma _{|A_i})$ then 
clearly, $U$ is definably isomorphic with a definable group of the form
$A_1\times \cdots \times A_{l-1}\times V\times A_{l+1}\times \cdots 
\times A_k$ where
$V$ the definable extension of $G$ by $A_l$ obtained from $c_l$. 
}
\end{fact}

\begin{lem}\label{lem limit c}
Let $(A,\gamma )$ be a definable $G$-module. 
Suppose that $G$ is a  
one-dimensional torsion-free definably connected
definable group, and let $c\in Z^n_{\N}$$(G,A,\gamma )$ (where $n>0$). If 
\[ \forall g_1,\dots ,g_{n-1}\in G,\,\lim_{k\into +\infty}c(g_1,\dots
,g_{n-1},k) \in A \]
then $c\in B^n_{\N}$$(G,A,\gamma ).$
\end{lem}

\pf
For each $g_1,\dots ,g_{n-1}\in G$ let 
\[ b(g_1,\dots ,g_{n-1}):=\lim_{k\into +\infty}c(g_1,\dots
,g_{n-1},k) \in A .\]
We have,
\[0=\gamma (g_1)(c(g_2,\dots ,g_{n+1})) +
\sum_{i=1}^n(-1)^ic(g_1,\dots,g_ig_{i+1},\dots ,g_{n+1})\, + \]
\[  (-1)^{n+1}c(g_1,\dots ,g_n) .\]
Taking the limit as $g_{n+1}\into +\infty$, we obtain (note that,
since $G$ is an ordered group $g_ng_{n+1}\into +\infty $ as
$g_{n+1}\into +\infty $)

\[ (-1)^nc(g_1,\dots ,g_n)=\gamma (g_1)(b(g_2,\dots ,g_n))\, + \]
\[ \sum_{i=1}^{n-1}(-1)^ib(g_1,\dots,g_ig_{i+1},\dots ,g_n)
+ (-1)^nb(g_1,\dots ,g_{n-1}) .\]
Therefore, $c$ is the coboundary of $(-1)^nb.$
\qed

\begin{lem}\label{lem splitting compact from one-dimensional}
Let $A$ be an abelian definably compact definable group 
such that $(A,\gamma )$ is a definable 
$G$-module. Suppose that $G$ is a one-dimensional torsion-free 
definably connected
definable group, then the action of $G$ on $A$ is trivial and
$Ext_{\N}$$(G,A,\gamma )$ is trivial. 
\end{lem}

\pf
This follows from lemma \ref{lem limit c} and 
the fact that $A$ is definably compact.
\qed

\begin{nrmk}\label{nrmk V}
{\em
Suppose that we have (definable) extensions 
$1\rightarrow A\rightarrow U\stackrel{\pi }{\rightarrow}
G\rightarrow 1 $ and $B\trianglelefteq G$ is definable then $C:=\pi
^{-1}(B)\trianglelefteq U$ and $A\trianglelefteq C$. 

Moreover, if we have a (definable) extension
$1\rightarrow B\rightarrow G\stackrel{j}{\rightarrow}
H\rightarrow 1$. Then we have (definable) extensions
$1\rightarrow C\rightarrow U\stackrel{j\pi }{\rightarrow}
H\rightarrow 1$ and 
$1\rightarrow A\rightarrow C\stackrel{\pi _{|C}}{\rightarrow}
B\rightarrow 1 $.
}
\end{nrmk}

\begin{lem}\label{lem splitting orthogonals}
Let $G$ be a one-dimensional definably connected torsion-free
definable group. For each $i\in \{1,\dots ,l\}$ let $A_i$ be a
definable group such that there are definable subgroups
$1=A^0_i<A^1_i<\cdots <A^{n_i}_i=A_i$ such that for each 
$j\in \{1,\dots ,n_i\}$, 
$A^j_i/A^{j-1}_i$ is definably isomorphic with a 
one-dimensional definably connected
torsion-free definable group with domain
$I_i$. Suppose that for each $i\in
\{1,\dots ,l\}$ there is no definable bijection between $G$ and $I_i$.
If $A:=A_1\times \cdots \times A_l$ and $(A,\gamma )$ is a definable $G$-module
then the action of $G$ on $A$ is trivial
and $Ext_{\N}(G,A,\gamma )$ is trivial.
\end{lem}

\pf
Its easy to see that the action of $G$ on $A$ is trivial. We now prove 
the rest. We have $Ext_{\N}$$(G,A,\gamma )$
$=H^2_{\N}$$(G,A,\gamma )=$
$H^2_{\N}$$(G,A_1,\gamma _{|A_1})\times \cdots \times 
H^2_{\N}$$(G,A_l,\gamma _{|A_l})$. We now use lemma \ref{lem
limit c} to conclude: let $c^i=(c^i_1,\dots ,c^i_{n_i})\in Z^2_{\N}$
$(G,A_i,\gamma _{|A_i})$ we 
will show that for each $j\in \{1,\dots ,n_i\}$ and each 
$g\in G$ the definable 
function $c^i_j(g,-):G\into A_i$ is such that
$\lim _{x\into +\infty }c^i_j(g,x)$ exists.
 
By the monotonicity
theorem $c^i_j(g,-)$ determines a definable bijection between an unbounded
interval in $G$ and an interval in $I_i^{n_i}.$ If this last interval is
bounded in $I_i^{n_i}$, then we are done. So suppose its unbounded.
But since we have
definable group structures on $I_i$ and on $G,$ this definable bijection
can be extended to a definable bijection between $I_i$ and $G,$ which is
a contradiction.
\qed

\begin{lem}\label{lem splitting lb}
Let $\I$ be an o-minimal structure,
$A_1=\cdots =A_l=(I,0,+)$ and $G=(I,0,\oplus )$  $\I$-definably connected 
one-dimensional torsion-free 
$\I$-definable groups. Let $A:=A_1\times \cdots \times A_l$ and
suppose that $(A,\gamma )$ is
an $\I$-definable $G$-module. If $\I$ is linearly bounded with 
respect to $+$ then the action of $G$ on $A$ is trivial and 
$Ext_{\N}$$(G,A,\gamma )$ is trivial.
\end{lem}

\pf
The fact that the action is trivial follows from the fact that
$\I$ is linearly bounded with respect to $+$. 
We now need to show that each $H^2_{\N}$$(G,A_i,\gamma _{|A_i})$ is
trivial. So we may assume
without loss of generality that $l=1$.
For this we use lemma \ref{lem limit c}. 
Let $c\in Z^2_{\N}$$(G,A,\gamma )$ be the definable
canonical 2-cocycle. 
Since $\I$ is linearly bounded with respect to $+$, there are 
$r_1,\dots ,r_l\in \Lambda (\I )$ such that for each $x,y\in G$ we have 
$c(x,y)=r_xy+o(x,y)$ where $r_x\in \{r_1,\dots ,r_l\}$ and $o:G\times
G\into A$ is a definable function such that for each $x\in G$ the
function $o_x:G\into A$, $y\into o(x,y)$ is bounded (in particular,
$\lim _{y\into +\infty}o(x,y)\in A$).

Let $g,h,k\in G$, and suppose $h$ is large enough so that
$r_h=r_{g\oplus h}=r$. Then by equation \ref{eq c(g,h) in Z2} we have 
$[r_g(h\oplus k)+o(g,h\oplus k)]-[r_gh+o(g,h)]+[o(h,k)-o(g\oplus h,k)]=0$.

And therefore $\forall g\in G,\, r_g=0$, since the above equality
implies that $r_g$ is bounded (take $k\into +\infty $). And so, 
$\forall g\in I,\,\lim_{h\into +\infty}c(g,h)\in I.$  
\qed

Recall the following important result.

\begin{fact}\label{fact ps}
{\em
\cite{ps}
Let $G$ be a definable group which is not definably compact, and let 
$\sigma :(a,b)\subseteq N\into G$ be a definable curve which is not 
completable in $G$ suppose without loss of generality that $\lim
_{x\into b^{-}}\sigma (x)$ does not exist in $G$. Let $I:=\sigma
((a,b))$. Then there is an induce order $<$ on $I$. Let $\M$ be an
$|N|^{+}$-saturated extension of $\N$; let $I^{\infty }=\{x\in
I^{\M}$$:\forall b\in I^{\N}$$x>b\}$ and for each $\alpha \in G^{\N}$
let $\mathcal{V}_{\alpha }$ be the infinitesimal neighbourhood of
$\alpha $ in $\M$, i.e., the intersection of all $\N$-definable
$V\subseteq G^{\M}$ of $\alpha $.

Define an equivalence relation on $G^{\N}$ by $\alpha T_{I}\beta \iff 
\mathcal{V}_{\alpha }*$$I^{\infty }=\mathcal{V}_{\beta }*$$I^{\infty }$
where $*$ is the group operation on $G$. Then the $T_{I}$-equivalence
class of the identity element of $G$ is a one-dimensional torsion-free
ordered definable subgroup $H_{I}$ of $G$ and the $T_{I}$-equivalence
classes are exactly the left cosets of $H$. 
}
\end{fact}

A corollary of the proof of fact \ref{fact ps} is the following remark which
shows the limitations of the method of fact \ref{fact ps} for 
finding one-dimensional torsion-free  ordered definable group. 
(We will use the notation of fact \ref{fact ps}).

\begin{nrmk}\label{nrmk ps}
{\em
Suppose that we have a definable extension  $1\rightarrow A\rightarrow
U\rightarrow G\rightarrow 1$ where $G$ is a one-dimensional
torsion-free ordered
definable group. Let $c\in Z^2_{\N}$$(G,A,\gamma )$ be the corresponding
definable $2$-cocycle. 

We know that we can assume that $U$ is a definable group with
domain $A\times G$ and with group operation given by
$(a,x)(b,y)=(a+b^g+c(x,y),xy)$. Suppose that
$dimU=n$. Then there is a definable open neighbourhood  of the identity
element of $U$ which is definably homeomorphic to a  definable open
subset $O\subseteq N^n$. Without loss of generality, we may assume
that $O\subseteq U$. For each $t\in G^{>1}$ let $B_t$ be an open
rectangular box such that $B_t\cap \{(0,x):x\in G\}=\{(0,x):t^{-1}<x<t\}$.
Let $t_0\in G^{>1}$ be such that $B_{t_0}\subseteq O$. For $1<t<t_0$
let $\bar{B_t}$ be the topological closure of $B_t$ in $O$
and let $bd(B_t)$ be its boundary.

For each $u\in G$ let $S_u:=\{(0,x)(0,u)^{-1}:x\in G^{\geq u}\}$. By
o-minimality, its easy to see that for all $1<t<t_0$, $S_u\cap
bd(B_t)\neq \emptyset $. Consider the following definable functions
$g:(1,t_0)\times G\into G$, $f:(1,t_0)\times G\into U$ and
$h:(1,t_0)\into U$ given by
\[
\begin{array}{clrc}
g(t,u):=\inf \{x\in G^{\geq u}:(0,x)(0,u)^{-1}\in bd(B_t)\},\\
f(t,u):=(0,g(t,u))(0,u)^{-1}\in bd(B_t)\,\,\,\,\,\, and\\
h(t):=\lim _{u\into +\infty }f(t,u).
\end{array}
\]

Note that, a simple calculation shows that
\[
\begin{array}{clrc}
f(t,u)=(-c(u,u^{-1})^{u^{-1}}+c(g(t,u),u^{-1})^{u^{-1}},g(t,u)u^{-1})
\end{array}
\]
By the proof of fact \ref{fact ps} (see claim 3.8.2 in \cite{ps}),
$Imh$ is a one-dimensional subset of $H_I$ where $I:=\{(0,x):x\in
G^{>1}\}$. Therefore, $U=A\rtimes _{\gamma }G$ iff $Imh\cap A\neq \{1\}$
iff for some $t\in (1,t_0)$, $\lim _{u\into +\infty }g(t,u)u^{-1}\neq
1.$ 
}
\end{nrmk}

Let $U$ be a definable abelian group 
of dimension two and with no definably compact parts. 
Lemma \ref{lem splitting orthogonals}
and lemma \ref{lem splitting lb} above show that $U$ is definably
isomorphic to a direct
product of two one-dimensional torsion-free
definable groups except possibly in the case where $U$ is 
a group definable in a
definable o-minimal expansion $\I$ of a real closed field
$(I,0,1,+,\cdot ,<)$ and we have an $\I$-definable extension
$1\rightarrow A\rightarrow U\rightarrow G\rightarrow 1$ where
$A=(I,0,+,<)$ and $G=(I,0,\oplus ,<)$ is a one-dimensional torsion-free
$\I$-definable group. 

\begin{lem}\label{lem ps problem}
Let $\I$ be an expansion of a real closed field. Suppose that we have
an $\I$-definable abelian extension  
$1\rightarrow A\rightarrow U\rightarrow G\rightarrow 1$
 where $A=(I,0,+,<)$ and $G=(I,0,\oplus ,<)$ is a one-dimensional torsion-free
$\I$-definable group. Then there is a $2$-cocycle 
$c\in Z^2_{\I}$$(G,A)$ be the corresponding to this
$\I$-definable  extension such that $U$ is $\I$-definably isomorphic to
$A\times G$ iff there is an $\I$-definable function $\alpha :G\into A$ 
such that
\[
\begin{array}{clrc}
\forall s\in G,\, \alpha
'(s)\frac{\partial \oplus }{\partial x}(0,s)=\alpha '(0)+\frac{\partial
c}{\partial x}(0,s).
\end{array}
\]
\end{lem}

\pf
Let
$t:G\into U$ be an $\I$-definable section, then by o-minimality there are
$g_0>\epsilon >0$ such that $t$ is $C^m$ on 
$(g_0\ominus \epsilon ,+\infty )$. Let
$s:G\into U$ be the $\I$-definable section given by: for all $g\in G$,
if $g>\ominus \epsilon $ then $s(g):=t(g\oplus g_0)t(g_0)^{-1}$ and if 
$g\leq \ominus \epsilon $ then 
$s(g):=s(\ominus g)^{-1}$. Then $s(0)=(0,0)$ and $s$ is $C^m$ on $G\setminus
\{\ominus \epsilon \}$. Let 
$c(g,h):=s(g)s(h)s(g\oplus h)^{-1}$ be the corresponding $\I$-definable 
$2$-cocycle, then $c$ is $C^m$ everywhere except
possibly on 
$\{\ominus \epsilon \}\times G \cup G\times \{\ominus \epsilon \}$. 

By 
fact \ref{fact extension 2-co-cycle}, $U$ is $\I$-definably isomorphic with 
$A\times G$ if and only if there is an $\I$-definable function 
$\alpha :G\into A$ with $\alpha (0)=0$ such
that the definable 
function $\beta :G\into U$, $\beta (s):=(\alpha (s),s)$ is a definable 
homomorphism, equivalently if and only if the $\I$-definable function 
$\alpha :G\into A$ satisfies
\[
\begin{array}{clcr}
\forall t,s\in G,\, \alpha (t\oplus s)=\alpha (t)+\alpha (s) +c(t,s)
\end{array}
\]
if and only if (to see this use also the fact that $U$ is abelian)
\[
\begin{array}{clcr}
\forall t,s\in G, \, \alpha '(t\oplus s)\frac{\partial \oplus
}{\partial x}(t,s)=\alpha '(t)
+ \frac{\partial c}{\partial x}(t,s).
\end{array}
\]

Putting $t=0$ in the second equation we get

\[
\begin{array}{clcr}
\forall s\in G, \, \alpha '(s)\frac{\partial \oplus }{\partial x}(0,s)
=\alpha '(0)+\frac{\partial c}{\partial x}(0,s).
\end{array}
\]

We show that this last equation is equivalent to the second equation:
putting $t\oplus s$ in the third equation we get $\alpha '(t\oplus s)
\frac{\partial \oplus }{\partial x}(0,t\oplus s)=\alpha '(0) +
\frac{\partial c}{\partial x}(0,t\oplus s);$ the associativity of
$\oplus $ implies that $\frac{\partial \oplus }{\partial x}(0,t\oplus
s)=\frac{\partial \oplus }{\partial x}(t,s)\frac{\partial \oplus
}{\partial x}(0,t);$ and since $c$ is a $2$-cocyle we get
$-\frac{\partial c}{\partial x}(t,s)\frac{\partial \oplus }{\partial
x}(0,t)+ \frac{\partial c}{\partial x}(0,t\oplus s)-\frac{\partial
c}{\partial x}(0,t)=0.$ From these equations together with the third
equation we get the second equation.
\qed

Using lemma \ref{lem ps problem} and results from \cite{sp}
and \cite{pss} we get:

\begin{cor}\label{cor splitting on the reals}
Let $\tilde{\RR}$ be an o-minimal expansion of $(\RR$$,0,+,<)$ the
additive group of real numbers. Then there is an o-minimal expansion 
$\hat{\RR}$ of $\tilde{\RR}$ such that every $\tilde{\RR}$-definable
abelian group with no $\tilde{\RR}$-definable
compact parts is $\hat{\RR}$-definably
isomorphic to a product of one-dimensional groups
$\hat{\RR}$-definably isomorphic to $(\RR$$,0,+)$ and
$(\RR$$^{>0},1,\cdot )$.
\end{cor}

Since the theory of (ordered) real closed fields has quantifier elimination in
the language of ordered rings, we have:

\begin{cor}\label{cor splitting on rcf}
Let $\mathbf{R}$=$(R,0,1,+,\cdot ,<)$ be a real closed field. 
Then every $\mathbf{R}$-definable abelian group with no
$\mathbf{R}$-definably compact parts
is $\mathbf{R}$-definably isomorphic to a product of one-dimensional groups
$\mathbf{R}$-definably isomorphic to $(R,0,+)$ and
$(R^{>0},1,\cdot )$.
\end{cor}

\end{subsection}

\begin{subsection}{The main theorem}
\label{subsection the main theorem}

We are now read to prove our main theorem.

\begin{thm}\label{thm the main theorem} 
Let $U$ be a definable solvable group. Then $U$ has a definable normal
subgroup $V$ such that $U/V$ is a definably compact definable solvable
group and 
$V=K\times W_1\times \cdots \times W_s
\times V'_1\times V_1 \times \cdots \times V'_k\times V_k$ 
where $K$ is the definably connected
definably compact normal subgroup of $U$ of maximal dimension and 
for each $j\in \{1,\dots ,s\}$ (resp., $i\in \{1,\dots ,k\}$) 
there is a semi-bounded o-minimal expansion $\J$$_j$ of a group 
(resp., an o-minimal expansion $\I$$_j$ of a real closed field)
definable in $\N$ such that there is no definable bijection between
a distinct pair among the $J_j$'s and $I_i$'s,
$W_j$ is a direct product of copies of the additive group of $\J$$_j$, 
$V'_i$ is a direct product of copies of the linearly bounded 
one-dimensional torsion-free $\I$$_i$-definable group and 
$V_i$ is an $\I$$_i$-definable group such that
$Z(V_i)$ has an $\I$$_i$-definable subgroup $Z_i$ such that
$Z(V_i)/Z_i$ is a direct product of copies of the linearly bounded
one-dimensional torsion-free $\I$$_i$-definable group and there are
$\I$$_i$-definable subgroups $1=Z^0_i<Z^1_i<\cdots <Z^{m_i}_i=Z_i$
such that for each $l\in \{1,\dots ,m_i\}$, $Z^l_i/Z^{l-1}_i$ is the
additive group of $\I$$_i$, and $V_i/Z(V_i)$ $\I$$_i$-definably embeds
into $GL(n_i,I_i)$.
\end{thm}

\pf
We prove this by induction on dimension of $U$. The result is clearly
true for dimension one.
So let $U$ be as above
and suppose that the result is true for solvable definable groups of 
lower dimensions than that of $U$.
 
Let $K$ be the definably compact, definably connected, definable
normal subgroup of $U$ of maximal dimension. This exists: let 
$K_1$ be a definably compact, definably connected, definable
normal subgroup of $U$ and let $U_1=U/K_1$. Let $K_2$ be a definably
compact, definably connected, definable
normal subgroup of $U_1$. Now apply remark \ref{nrmk V} and 
let $K_3$ be the definable normal
subgroup of $U$ which is a definable extension of $K_2$ by $K_1$.
$K_3$ is a definably compact, definably connected, definable
normal subgroup of $U$ with $dimK_3\geq dim K_1$. Repeating this
process finitely many times we obtain $K$.

Let
$U':=U/K$. Then $U'$ is definable and has a definable normal (solvable)
subgroup with no definably compact parts, for otherwise the only definable
normal (solvable) subgroups of $U'$ would be definably compact and so
by remark \ref{nrmk V}, $K$ would not be maximal. Let $Y$ be the
maximal definable normal subgroup of $U'$ with no definably compact
parts (this exists by an argument similar to that above). Then
$U'/Y$ is definable and  definably compact, for otherwise by the
induction hypothesis $U'/Y$ would have a definable normal subgroup
with no definably compact parts and by remark \ref{nrmk V}, $Y$ would not be
maximal.
Now apply remark \ref{nrmk V} and let $V$ be the definable normal
subgroup of $U$ which is a definable extension of $Y$ by $K$.
Note that $U/V=U'/Y$ and so, it is definably compact. 
 
We now proceed with the proof, we will same times use fact
\ref{fact Al split}. Its easy to verify that each time we do this,
all the hypothesis are satisfied.
By repeated application of remark \ref{nrmk V}, lemma \ref{lem 
splitting compact from one-dimensional}
and fact \ref{fact Al split} we see that $V=K\times Y$. By
induction hypothesis and by repeated application of 
remark \ref{nrmk V}, lemma \ref{lem splitting orthogonals} 
and fact \ref{fact Al split}
we see that $Y=Y_1\times \cdots
\times Y_r$ where for each $i\in \{1,\dots ,r\}$ 
there are definable subgroups
$1=Y^0_i<Y^1_i<\cdots <Y^{n_i}_i=Y_i$ such that for each 
$j\in \{1,\dots ,n_i\}$, 
$Y^j_i/Y^{j-1}_i$ is definably isomorphic with a 
one-dimensional definably connected
torsion-free definable group with domain
$I_i$ and for $j\neq i$,
there is no definable
bijection between $I_i$ and $I_j$. For each $i\in \{1,\dots ,r\}$ let
$\I$$_i$ be the definable structure induced by $\N$ on $I_i$.

If $\I$$_i$ is a semi-bounded o-minimal expansion of a group then we make
$\J$$_i$:=$\I$$_i$, and $W_i:=Y_i$. And by induction on $dim W_i$ and 
applying (if needed) remark \ref{nrmk V} and lemma \ref{lem splitting
lb} several times we are
done. So assume that
$\I$$_i$ is an o-minimal expansion of a real closed field. Let $Y'_i$
be the maximal $\I$$_i$-definable normal subgroup of $Y_i$ which is a
direct product of copies of the linearly bounded one-dimensional
torsion-free $\I$$_i$-definable group. Then $X_i:=Y_i/Y'_i$ is definable
and by repeated application of remark \ref{nrmk V}, 
lemma \ref{lem splitting lb} and fact \ref{fact Al split}
we see that $Y_i=Y'_i\times X_i$. Now put $V'_i:=Y'_i$ and $V_i:=X_i$.
The fact that $Z(V_i)$ is as described is proved in the same way. The
fact that $V_i/Z(V_i)$ $\I$$_i$-definably embeds into some
$GL(n_i,I_i)$ is proved in \cite{opp}.
\qed

Corollary \ref{cor minimal reduct} 
below is an adaption of an argument due to Iwasawa
(see the proof of lemma 3.4 \cite{i}). We will need the following
result from \cite{s}. Recall that a definable group $G$ is 
{\it monogenic} if there is $g\in G$ such that the smallest definable
group containing $g$ (which exists by DCC) is $G$.

\begin{fact}\label{fact monogenic quotient}
{\em \cite{s}
Let $A\trianglelefteq U$ be definable groups. If $A\subseteq Z(U)$ and
$U/A$ is monogenic then $U$ is abelian.
}
\end{fact}

\begin{cor}\label{cor minimal reduct}
Let $\I$$=(I,0,1,+,\cdot ,<,\dots )$ be an o-minimal expansion of a
real closed field 
and suppose that there are no Peterzil-Steinhorn $\I$-definable
groups. Let $U$ be an $\I$-definable solvable group with no $\I$-definable 
compact parts. Then $U$ is $\I$-definably
isomorphic to a group definable of the form 
$U'\times G_1\cdots G_k\cdot G_{k+1}\cdots G_l$ where $U'$ is a direct product
of copies of the linearly bounded one-dimensional torsion-free
$\I$-definable group, for each $i\in \{1,\dots ,k\}$, $G_i=(I,0,+)$
and for each $i\in \{k+1,\dots ,l\}$, $G_i=(I^{>0},1,\cdot )$. In
particular, $G:=G_1\cdots G_k\cdot G_{k+1}\cdots G_l$ $\I$-definably
embeds into some $GL(n,I)$ and $U$ is $\I$-definably isomorphic to a
group definable in one of the following  reducts $(I,0,1,+,\cdot
,\oplus )$,  $(I,0,1,+,\cdot ,\oplus ,e^t)$ or 
$(I,0,1,+,\cdot ,\oplus , t^{b_1},\dots
,$$t^{b_r})$ of $\I$ where $(I,0,\oplus )$ is the Miller-Starchenko
group of $\I$, $e^{t}$ is the $\I$-definable exponential map (if it
exists), and the $t^{b_j}$'s are $\I$-definable power
functions. Moreover, if $U$ is nilpotent then $U$ is $\I$-definably 
isomorphic to a
group definable in the reduct $(I,0,1,+,\cdot ,\oplus )$ of $\I$. 
\end{cor}

\pf
By theorem \ref{thm the main theorem}, we may assume that $U=U'\times
G$ where $U'$ is the maximal $\I$-definable normal subgroup of $U$
which is a product of copies of 
the linearly bounded one-dimensional torsion-free
$\I$-definable group and $G$ is as described there. Furthemore, since
there are no Peterzil-Steinhorn $\I$-definable groups, every 
$\I$-definable abelian group with no $\I$-definably compact parts is a
direct product of one-dimensional torsion-free
$\I$-definable groups and therefore by an argument similar to those
used in the proof of theorem \ref{thm the main theorem} we can assume
that $Z(G)$ is a direct product of copies of additive group of $\I$.
Moreover, by an argument similar to that used in the proof of theorem 
\ref{thm the main theorem} (substitute ``$\I$-definably compact
$\I$-definable group'' by ``linearly bounded one-dimensional torsion-free
$\I$-definable group''), there are $\I$-definable subgroups 
$1=H_0\trianglelefteq H_1\trianglelefteq
\cdots \trianglelefteq H_{n+1}=G$ such that for each $i\in \{ 1, \dots
,n\}$, $H_i$ is the smallest definable normal subgroup of $H_{i+1}$
such that $H_{i+1}/H_i$ is abelian, 
$H_i/H_{i-1}$ is a direct product of copies  of additive group of $\I$
and $H_{n+1}/H_n$ is a direct product of copies (possibly zero
copies) of the linearly bounded
one-dimensional torsion-free
$\I$-definable group. 

Let $\bar{G}:=G/Z(G)$. Since $\bar{G}$ $\I$-definably embeds into some
$GL(k,I)$, by \cite{pps3} and the remark above,
$\bar{G}=\bar{G}_1\cdots \bar{G}_{\bar{k}}\cdot \bar{G}_{\bar{k}+1}\cdots
\bar{G}_{\bar{l}}$ where for each $i\in \{1,\dots ,\bar{k}\}$, 
$\bar{G}_i=(I,0,+)$
and for each $i\in \{\bar{k}+1,\dots ,\bar{l}\}$,
$\bar{G}_i=(I^{>0},1,\cdot )$.

Let $N$ be the $\I$-definable extension of 
$\bar{G}_1\cdots \bar{G}_{\bar{k}}\cdot \bar{G}_{\bar{k}+1}\cdots
\bar{G}_{\bar{l}-1}$ by $Z(U)$ (and therefore $G/N$ is a one dimensional 
torsion-free $\I$-definably connected
$\I$-definable group).
By induction its enough to show that $G$ contains an 
$\I$-definable subgroup 
$H$ ($\I$-definably isomorphic with $G/N$) such that $G=NH$ and 
$H\cap N=1$.

We prove this by induction on $\bar{l}$. Note that if $\bar{l}=0$ 
or $\bar{l}=1$, then $G$ is abelian (in the second case by fact  
\ref{fact monogenic quotient}) and so the claim holds by assumption.
Assume that the claim is true all $\I$-definable groups
with no $\I$-definably compact parts and with lower $\bar{l}$.

Suppose that $N$ contains a proper 
$\I$-definable normal subgroup $N_1$ of $G$. 
By induction applied to $G/N_1$
there is an $\I$-definable subgroup $G_1$ such that $G=NG_1$, 
$G_1\cap N=N_1$
and $G_1/N_1=G/N.$ Again the induction assumption for $G_1$ and $N_1$ 
gives us an  $\I$-definable subgroup $H$ such that $G_1=N_1H$ and 
$H\cap N_1=1$. This $H$ satisfies the claim.

We can therefore assume that $N$ has no proper $\I$-definable  
subgroup which is normal in  $G$. 
If $N$ is in the centre of $G$ then by fact  
\ref{fact monogenic quotient} $G$ is abelian and by assumption 
the claim is proved. If $N$ is not in the centre 
of $G$ then, using the decomposition series $1=K_0\trianglelefteq K_1
\trianglelefteq \cdots \trianglelefteq K_{m+1}=N$ of $N$ like the one
we got above for $G$, we see that
$N$ must be a direct product of $k$ copies of the additive group of $\I$. 
$N$ is therefore an $\I$-definable
$G$-module under conjugation and we have a natural $\I$-definable
homomorphism $A:G\into GL(k,I)$. 
$G/N$ $\I$-definably embeds in $GL(k,I)$. We show that
that there is $g\in G$ such that $det(A(g)-Id)\neq 0$ and so
$[N,g]=N$. Since $N$ is not in the centre of $G$, there is $g\in G$
which does not commute with some element in $N$. Let $N'$ be the
eigen-space for the value $1$ of the matrix $A(g)$. Since $A(G)$ is
abelian, $N'$ is invariant under all the $A(h)$. But this means that
the $\I$-definable subgroup $N'$ of $N$ is normal in $G$ and therefore
by the assumption we must have either $N'=N$ or $N'=1$. The first case
does not hold since $g$ does not commute with some element of
$N$. Therefore $N'=1$, $det(A(g)-Id)\neq 0$ and $[N,g]=N$.

Now take an arbitrary element $y\in G$ and put $z:=gyg^{-1}y^{-1}$.
Since $G/N$ is abelian, we have $z\in N$. Take $u\in N$ such that 
$z=gug^{-1}u^{-1}$ and put $v:=u^{-1}y$. It follows that 
$gv=vg$ and so $G=NC_G(g)$. If $x\in C_G(g)\cap N$, then
$gxg^{-1}x^{-1}=1$ and $det(A(g)-Id)\neq 0$ implies that $x=1$ , i.e., 
$C_G(g)\cap N=1$.

We have $G:=G_1\cdots G_k\cdot G_{k+1}\cdots G_l$. An induction on $l$
shows that $G$ $\I$-definably
embeds into some $GL(n,I)$ and $G$ is $\I$-definably isomorphic to a
group definable in one of the following  reducts $(I,0,1,+,\cdot )$,  
$(I,0,1,+,\cdot ,e^t)$ or 
$(I,0,1,+,\cdot , t^{b_1},\dots
,$$t^{b_r})$ of $\I$ where $e^{t}$ is the $\I$-definable exponential map (if it
exists), and the $t^{b_j}$'s are $\I$-definable power
functions. If $U$ is nilpotent then $G$ is nilpotent and by 
\cite{pps3}, $G$ is $\I$-definably 
isomorphic to a group definable in the reduct $(I,0,1,+,\cdot )$ of $\I$.
\qed

\begin{nrmk}\label{nrmk minimal reduct}
{\em
\cite{pps3}
There are solvable linear
groups $U$ and $V$ definable in o-minimal expansions of
$(\RR$$,0,1,+,\cdot ,<)$
by the $exp $ and $t^r$ respectively, 
such that $U$ (resp., $V$) is not isomorphic (even abstractly) to a 
definable in o-minimal expansions of $(\RR$$,0,1,+,\cdot ,<)$
by some $t^s$
(resp., real semialgebraic group): Let $A=(\RR$$^2,0,+)$, $G=(\RR$$,0,+)$ and
$H=(\RR$$^{>0},1,\cdot )$. Let $U=A\rtimes _{\alpha }G$ and
$V=A\rtimes _{\beta }H$, where $\alpha (t)(a,b)=(exp(t)a+texp(t)b,exp(t)b)$ and
$\beta (t)(a,b)=(ta,t^rb)$. 
}
\end{nrmk}

We end this subsection with the following result from \cite{ps} which
shows that definable abelian groups are not necessarily the direct
product of a definable abelian group with no definably compact parts and a
definably compact definable abelian group.

\begin{fact}\label{fact tmn}
{\em
\cite{ps}
Let $\tilde{\RR}$$:=(\RR$$,0,1,+,<)$. Then
for $m,n\in \NN$ and $L$ an integral lattice in $\RR$$^n$
there are $\tilde{\RR}$-definable abelian
groups $T(m,n,L)$ and $T(n,L)$ 
with dimensions $m+n$ and $n$ respectively, such 
that we have an $\tilde{\RR}$-definable extension
$1\rightarrow (\RR$$^m,0,+)\rightarrow T(m,n,L)\rightarrow T(n,L)
\rightarrow 1$. Moreover, if $L$ is "generic" then $(\RR$$^m,0,+)$ 
does not have an $\tilde{\RR}$-definable complement in $T(m,n,L)$ and
$T(n,L)$ does not have $\tilde{\RR}$-definable 
infinite proper subgroups. 

The same result holds in $(\RR$$,0,1,+,\cdot ,<)$.
}
\end{fact}

\end{subsection}
\end{section}

\begin{section}{The Lie-Kolchin-Mal'cev theorem}
\label{section lie kolchin malcev theorem}

\begin{subsection}{More on definable $G$-modules}
\label{subsection more on definable g modules}

In this subsection we will describe definable $G$-modules,
generalising a result from \cite{mmt} describing faithful irreducible
definable $G$-modules.

\medskip
\textbf{Notation:}
Let $(A,\gamma )$ be a definable $G$-module, for $i=1,\dots ,m$
let $(A_i,\gamma _i)$ be a definable $G_i$-module
and let $(B,\gamma )$ be a definable trivial
$G$-submodule. We write $(G,A,\gamma )=
(G_1,A_1,\gamma _1)\times $
$\cdots \times (G_m,A_m,\gamma _m)$ if 
$G=G_1\times \cdots \times G_m$, $A=B\times A_1\times \cdots \times A_m$ 
and for all $g=(g_1,\dots ,g_m)\in G$, for all
$a=(b,a_1,\dots ,a_m)\in A$ we have 
$\gamma (g)(a)=(b,\gamma _1(g_1)(a_1),\dots ,\gamma _m(g_m)(a_m))$.
Recall also that $\bar{G}$ denotes $G/Ker\gamma $ and we have a natural
definable $\bar{G}$-module $(A,\bar{\gamma})$. Also, $\bar{A}:=A/A^G$
and we have a natural definable $G$-module 
$(\bar{A}, \gamma _{\bar{A}})$.

\begin{thm}\label{thm non trivial}
Let $(U,\gamma )$ be a definable non trivial $G$-module where, 
$U$ and $G$ are infinite definably connected definable groups. Then
there is a definable subgroup $V$ of $U$ of the form $K\times W\times
V_1\times \cdots \times V_m$ such that $(K,\gamma )$ is the
maximal definably
connected definably compact trivial $G$-submodule of $(U,\gamma )$,
$(W,\gamma )$ is the maximal trivial $G$-submodule of $(U,\gamma )$ with no
definably compact parts and for each $i\in \{1,\dots ,m\}$ there is a definable
o-minimal expansion $\I$$_i$ of a real closed field $I_i$ such that
if $j\neq i$ then there is no definable bijection between $I_i$
and $I_j$ and (i) $(V_i,\gamma )$ is an $\I$$_i$-definable non trivial 
$G$-submodule of $(U,\gamma )$ with no definably compact parts and non $\I$$_i$
-linearly bounded, (ii) $(\bar{G},V,\bar{\gamma })
=(G_1,V_1,\gamma _1)\times \cdots \times (G_m,V_m,\gamma _m)$ 
where for each $i$, $G_i$ is definably isomorphic to an
$\I$$_i$-definable subgroup of some $GL(k_i,I_i)$,
$(V_i,\gamma _i)$ is a faithful definable $G_i$-module,  
$\bar{V_i}=(I_i^{>_i0_i},1_i, \cdot _i)^{l_i}\times (I_i,0_i,+_i)^{n_i}$ 
and (iii) $(G_i,\bar{V_i},\gamma _{i|\bar{V_i}})=
(H_i^1,U_i^1,\alpha _i^1)\times \cdots \times (H_i^{m_i},
U_i^{m_i},\alpha _i^{m_i})$ where for each $j\in \{1,\dots ,m_i\}$,
$(U_i^j,\alpha _i^j)$ is a $I_i$-semialgebraic faithful and 
irreducible $H_i$-module and $H_i/Z(H_i)$ is a direct product of 
$I_i$-semialgebraic non abelian $I_i$-semialgebraically simple groups.
Moreover, $(U/V,\gamma _{|U/V})$ is a definably compact trivial 
definable $G$-module. 
\end{thm}

\pf
We will refer to the notation of theorem \ref{thm the main theorem}. 
Its clear from theorem \ref{thm the main theorem} the existence of $V$ with
$K$ and $W$ with the properties mentioned, so to finish the prove
of (i) its enough to show that there is a definable non linearly
bounded and with no definably compact parts non trivial $G$-submodule. 
Suppose this is not the
case. Then by theorem \ref{thm compact G modules} and by \cite{ms} it
follows that $V$ is
contained in $U^G$ and so $\bar{U}:=U/U^G$ is a definably compact, 
definably connected definable group. 
By theorem \ref{thm compact G modules} 
$(\bar{U},\gamma _{\bar{U}})$ is a trivial $G$-module and so 
$\forall g\in G\forall \bar{u}\in \bar{U},\, 
\gamma _{\bar{U}}(g)(\pi ^{-1}(\bar{u}))\subseteq 
\pi ^{-1}(\bar{u})$ (where, $\pi :U\into\bar{U}$ is the natural 
projection) and therefore, if $B$ is
an infinite minimal definable subgroup of $\bar{U}$ we have a
definable family $\Gamma :G\times B\into U^G$ of definable
homomorphisms from $B$ into $U^G$ given by, $\forall g\in G\forall
b\in B,\, \Gamma (g)(b):=\gamma _{\bar{U}}(g)(x)-x$ for some 
$x\in \pi^{-1}(b)$. Now, since
$(U,\gamma )$ is a non trivial definable $G$-module,
by an argument similar to that in the proof of
theorem \ref{thm compact G modules} we get a contradiction.

We now prove (ii).
Now let $k_i:=dimV_i$. By corollary 2.21 and fact 2.24 in \cite{pps1}
we have, after fixing a basis for the tangent space
of each $V_i$ a definable homomorphism 
$G\into GL(k_1,I_1)\times \cdots \times GL(k_m,I_m)$ 
given by 
$g\into (d_0(\gamma _{|V_1}(g)),\dots ,d_0(\gamma _{|V_m}(g))$ and with
kernel $Ker\gamma $. This shows that $\bar{G}=G_1\times \cdots \times
G_m$ where each $G_i$ is definably
isomorphic with an $\I$$_i$-definable subgroup of $GL(k_i,I_i)$.
Since $G$ is definably connected, each $G_i$ is infinite and since for
$j\neq i$ there is no definable bijection between $I_i$ and $I_j$,
we have $G_i\subseteq Ker\gamma _{|V_j}$, so 
to prove the first part of (ii), take $\gamma _i:=\bar{\gamma }_{|V_i}$.

Consider $G_i$ as an $\I$$_i$-definable group and consider
the $\I$$_i$-definable group $V_i\rtimes _{\gamma _i}G_i$
whose center is $V_i^{G_i}\times (Ker\gamma _i\cap Z(G_i))$
$=V_i^{G_i}\times \{1\}$. By \cite{opp}, when have that
$V_i\rtimes G_i/(V_i^{G_i}\times \{1\})$ is 
$\I$$_i$-definably 
isomorphic with an $\I$$_i$-definable subgroup of some $GL(l_i,I_i)$
and so by \cite{pps3} 
$\bar{V_i}=(I_i^{>_i0_i},1_i,\cdot _i)^{l_i}\times (I_i,0_i,+_i)^{n_i}$.

We will now prove (iii). We clearly have $\bar{V_i}=U_i^0\times
U_i^1\times \cdots \times U_i^{m_i}$ where 
$(U_i^0,\gamma _{i|\bar{V_i}})$ is 
trivial $\I$$_i$-definable $G_i$-submodule of 
$(\bar{V_i},\gamma _{i|\bar{V_i}})$, and for each 
$j\in \{1,\dots ,m_i\}$, $(U_i^j,\gamma _{i|\bar{V_i}})$ is a faithful
and irreducible $\I$$_i$-definable $G_i$-submodule of 
$(\bar{V_i},\gamma _{i|\bar{V_i}})$ and therefore each such $U_i^j$
is a vector space over the real closed field $I_i$. O-minimality
implies that the action of $G_i$ on $U_i^j$ is by vector space
automorphisms and so we can easily get $H_i^j$ and $\alpha _i^j$ 
satisfying the first part of (iii). The rest is proved in 
proposition 1.3 \cite{mmt}.  
\qed

Peterzil and Starchenko proved in \cite{pst2}, using the theory
of $\bigvee $-definable groups and assuming 
that $\N$ has definable Skolem functions, that if $\mathbf{U}$$:=(U,\cdot )$
is a definable group which is not abelian-by-finite, then a real
closed field is interpretable in $\mathbf{U}$. Here we get the following.

\begin{cor}\label{cor defining a field}
Let $U$ be a definable group which is not abelian-by-finite. Then a
real closed field is definable in $(N,<,U,\cdot )$.
\end{cor}

\pf
Suppose that $U$ is definably connected.
Let $R(U)$ be the maximal definably connected
definable normal solvable subgroup of $U$.
If $R(U)$ is abelian then it is a definable $U$-module under
conjugation and if it is non-trivial we can apply theorem \ref{thm 
non trivial}, otherwise we have $Z(U)=R(U)$ and
$U/Z(U)$ is an infinite definably semi-simple definable group
and the result follows from \cite{pps1} and \cite{pps2}.

So suppose that $R(U)$ is not abelian. Since it is solvable, it has
a definable abelian normal subgroup $X$ such that $Z(R(U))\leq
X$ and $X/Z(R(U))$ is an infinite definable abelian group. $X$
is a non-trivial definable $R(U)$-module and we can apply
theorem  \ref{thm non trivial}.
\qed

\end{subsection}

\begin{subsection}{The Lie-Kolchin-Mal'cev theorem}
\label{subsection the lie kolchin malcev theorem}

Let $G$ be a definable group and $X$ a subset of $G$. By DCC on
definable subgroups, the intersection of all definable subgroups of
$G$ containing $X$ is a definable subgroup of $G$. This is the
smallest definable subgroup of $G$ containing $X$ and we denote it by
$d(X)$ and call it the {\it definable subgroup of $G$ generated by
$X$}.

\begin{lem}\label{lem on d(X)}
Let $G$ be a definable group. Then the following holds: (1)
The operator $d$ is a closure operator i.e., for all subsets $X,Y$ of
$G$ we have
$X\subseteq d(X)$, if
$Y\subseteq X$ then $d(Y)\leq d(X)$ and $d(d(X))=d(X)$. 
(2) If the
elements of $X\subseteq G$ commute with each other, then $d(X)$ is abelian. 
(3) If a
subgroup $A\leq G$ normalises the subset $X\subseteq G$, 
then $d(A)$ normalises $d(X)$.
(4) If $X,Y\leq G$ then $[d(X),d(Y)]\leq d([X,Y])$
in particular, a subgroup $H\leq G$ is solvable (resp., nilpotent) of class
$n$ iff $d(H)$ is also solvable (resp., nilpotent) of class $n$.
\end{lem}

\pf
$(1)$ is trivial. For $(2)$ and $(3)$ see the proof of lemma 5.35 in
\cite{bn}. As for $(4)$, the proof in \cite{bn} for the finite Morley
rank analogue (see corollary 5.38 and lemma 5.37 in \cite{bn})
works in our case using the following result (which is a consequence
of DCC): if $G$ is a definable
group and, $H$ is a definable normal subgroup of $G$, $A$ is a subgroup of
$G$ containing $H$ and $Y$ is a subset of $G$ containing $H$ are such that 
$A/H=C_{G/H}(Y/H)$, then $A$ is definable. 
\qed

\begin{lem}\label{lem lie kolchin malcev}
Let $G$ be a definable group. (1) If $G$ is definably connected then,
every finite normal subgroup is contained in $Z(G)$ and if $Z(G)$ is
finite then $G/Z(G)$ is centerless. (2) If $G$ is infinite and
nilpotent then $Z(G)$ is infinite. (3) If $G$ is infinite solvable but
not nilpotent then $G$ has an infinite proper maximal
normal definable subgroup $H$ such that $G/H$ is abelian.
\end{lem}

\pf
(1) is the o-minimal
analogue of corollary 1 in \cite{n} and lemma 6.1 in \cite{bn}. The
proof is the same. (2) is the o-minimal
analogue of lemma 6.2 in \cite{bn} again the proof is the same. (3) is
proved by an argument
contained in the proof of theorem 2.12 in \cite{pps2}. 
\qed

We are now ready to prove the o-minimal version of the
Lie-Kolchin-Mal'cev theorem. The proof is a modification of that in
\cite{n} for the finite Morley rank case.

\begin{thm}\label{thm lie kolchin malcev}
If $U$ is a definably connected
definable solvable group, then $U^{(1)}$ is a $\bigvee
$-definable nilpotent normal subgroup and $d(U^{(1)})$ is a definable 
nilpotent normal subgroup. 
\end{thm}

\pf
Let $U$ be a minimal counter-example, so both $U^{(1)}$ and
$d(U^{(1)})$ are not nilpotent.

\medskip
{\it Claim (1):} We can assume that $Z(U)=Z(U^{(1)})=1$.

\medskip
{\it Proof of Claim (1):} 
The fact that we may assume $Z(U)=1$ follows from 
 $(U/Z(U))^{(1)}$ $=U^{(1)}Z(U)/Z(U)$
$\simeq U^{(1)}/U^{(1)}\cap Z(U)$ $\supseteq U^{(1)}/Z(U^{(1)})$ (because 
$U^{(1)}\cap Z(U)\leq Z(U^{(1)})$) and so $Z(U)$ is finite and we
can substitute $U$ by $U/Z(U)$ which is centerless by lemma \ref{lem
lie kolchin malcev}.

By lemma \ref{lem quotient by AG and ker are 
definable} $U/C_U(U^{(1)})$ is definable. We have: 
$(U/C_U(U^{(1)}))^{(1)}$
$=U^{(1)}C_U(U^{(1)})/C_U(U^{(1)})$ $\simeq U^{(1)}/U^{(1)}\cap
C_U(U^{(1)})$ $=U^{(1)}/Z(U^{(1)})$. And so, if $C_U(U^{(1)})$ is
infinite then 
$(U/C_U(U^{(1)}))^{(1)}$ is nilpotent and so $U^{(1)}$ is also
nilpotent. Therefore, $C_U(U^{(1)})$ is finite and by lemma \ref{lem
lie kolchin malcev} we have
$Z(U^{(1)})\subseteq  C_U(U^{(1)})\subseteq Z(U).$

\medskip
{\it Claim (2):} $U^{(1)}$ and $d(U^{(1)})$ are torsion-free.

\medskip
{\it Proof of Claim (2):} 
We have 
$U^{(1)}\leq d(U^{(1)})\leq W_1\times \cdots \times W_s\times
V'_1\times V_1\times \cdots \times V'_k\times V_k$ and this last group
is torsion-free (this can be proved by induction on dimension and
using equation (\ref{eq group of h})). 

\medskip
{\it Claim (3):} There is an infinite definable abelian normal subgroup $A$ of
$U$ which is an irreducible faithful definable $U/C_U(A)$-module
under conjugation.

\medskip
{\it Proof of Claim (3):} Since $U$ is not nilpotent, by lemma
\ref{lem lie kolchin malcev}
$U$ has an infinite proper maximal
normal definable subgroup $X$ such that $U/X$ is abelian. Therefore, 
$d(U^{(1)})$ is an infinite definable normal proper subgroup of $U$
and so $U^{(2)}\subseteq d(U^{(1)})^{(1)}\subseteq
d(d(U^{(1)})^{(1)})$ is nilpotent and infinite (for otherwise,
$U^{(2)}$ is finite and since by claim (2) $U^{(1)}$ is
torsion-free, 
$U^{(2)}=1$ and $U^{(1)}$ would be abelian). Now by lemma \ref{lem lie
kolchin malcev}, $Z(d(d(U^{(1)})^{(1)}))$ is
infinite. Now let $A$ be an infinite definable normal subgroup of $U$
contained in $Z(d(d(U^{(1)})^{(1)}))$ and minimal for these
properties. Note that we have $U^{(2)}\leq C_U(A)$ and $U/C_U(A)$ is
infinite because otherwise we would have $A\leq Z(U)=1$. By minimality
of $A$, $A$ is an irreducible faithful definable $U/C_U(A)$-module
under conjugation. 

\medskip
By theorem \ref{thm non trivial}, $U/C_U(A)$
is abelian (since is solvable) and therefore we have
$1=(U/C(U(A))^{(1)}=U^{(1)}C_U(A)/C_U(A)$ $\simeq
U^{(1)}/C_{U^{(1)}}(A)$ and therefore, $U^{(1)}=C_{U^{(1)}}(A)$
i.e., $A\leq Z(U^{(1)})=1$ contradicting claim (3).
\qed

We finish this subsection with the following result on definable
nilpotent groups. Recall that a group $G$ is the central product of
two subgroups $H$ and $K$ if $G=HK$, $H$ and $K$ are normal and $H\cap
K\leq Z(G)$. We denote this by $G=H*K$. $H$ is {\it
divisible} if for every $n\in \NN$ and
every $x\in H$ there is $y\in H$ such that $y^n=x$.

\begin{thm}\label{thm definable nilpotent}
Let $B$ be a definable nilpotent group. Then $B=B^0*F$ for some finite
subgroup $F$ and $B^0$ is divisible. 
Moreover, if $B$ is abelian then $B=B^0\times F$ and
if we have an extension
$1\rightarrow A\stackrel{i}{\rightarrow}U\stackrel{j}{\rightarrow}
G\rightarrow 1$ where $A$ is an abelian definable group and $G$ is a
finite group then $U$ is definable.
\end{thm}

\pf
We will first prove the second part of the theorem.
Let
$1\rightarrow A\stackrel{i}{\rightarrow}U\stackrel{j}{\rightarrow}
G\rightarrow 1$
be an extension, where $A$ is an abelian definably connected
group and $G$ is a finite group.
Its clear that every abelian definably connected definable group $H$ is 
divisible: for every $n\in \NN$, the kernel of the
homomorphism $H\into nH$, $h\into nh$ is a definable subgroup of $H$
with  bounded exponent, and therefore by \cite{s} 
is finite and so $nH=H$.
An argument similar to that of lemma \ref{lem limit c} where we use
$\sum_{k\in G}$ instead of $\lim _{g\into +\infty}$ show that
if $A$ is a definable abelian connected group (and therefore
divisible) and $G$ is a finite group, then $H^n(G,A)$ is trivial and
this proves the second part of the theorem.

Let $B$ be a minimal counterexample to the first part of the theorem.
Then by the above, $B$ is not abelian-by-finite, $Z(B)^0$ is infinite
and $B/Z(B)^0$ is infinite. And so $B/Z(B)^0=(B/Z(B)^0)^0*F$. Let
$H$ and $K$ be definable normal subgroups of $B$ such that $H/Z(B)^0$
$=(B/Z(B)^0)^0$ and $K/Z(B)^0=F$. We have $K\neq B$ and by induction
$K=K^0*F_1$. Now we have $B=(K^0H)*F_1$ and by exercise 14, page 6
\cite{bn}, $K^0H$ is divisible and therefore, also definably
connected, i.e., $K^0H=B^0$.
\qed

\end{subsection}
\end{section}

\begin{section}{Existence of strong definable choice}
\label{section existence of strong definable choice}

Here we finally prove that definable groups have strong definable choice.

\begin{thm}\label{thm definable choice}
Let $U$ be a definable group and let $\{T(x):x\in X\}$ be a definable
family of non empty definable subsets of $U$. Then there is a definable
function $t:X\into U$ such that for all $x,y\in X$ we have
$t(x)\in T(x)$ and if $T(x)=T(y)$ then $t(x)=t(y)$.
\end{thm}

\pf
Let $R(U)$ be the maximal definable solvable normal subgroup of $U$.
Then $U/R(U)$ is definable and by \cite{pps1} it 
has the property stated in the theorem. On the other hand, 
there is a definable section $s:U/R(U)\into U$ 
and so $U$ is definably isomorphic
to a definable group with domain $R(U)\times U/R(U)$ and so, it
is sufficient to prove the theorem for definable solvable groups.
By theorem \ref{thm the main theorem} and an argument similar to the
one above,  
the result is true for definable solvable groups if it is
true for definably compact definable abelian groups. 

So let $U$ be a definably compact definable group and let
$\{T(x):x\in X\}$ be a definable family of non empty definable subsets
of $U$.
First note that the (induced) topology for the definable family
$T=\{T(x):x\in X\}$ is uniformly definable. Let
$\bar{T}:=\{\bar{T(x)}:x\in X\}$ where $\bar{T(x)}$ is the closure
of $T(x)$ in $U$.
Suppose that $U\subseteq N^m$ and for each $l\in \{1,\dots ,m\}$
let $\pi _l:N^m\into N^l$ be the projection onto the first $l$
coordinates.
For each $x\in X$  let $Y_m(x):=\bar{T(x)}$  and for each
$i\in \{1,\dots ,m-1\}$ let
$Y_i(x):=\{a\in \pi _i(U): \{ (a,b)\in \pi _{i+1}(U)\}\cap Y_{i+1}(x)
\neq \emptyset \}$ and let $Y_i:=\bigcup _{x\in X}Y_i(x)$. 
Note that for each $x\in X$ and each 
$a\in Y_{m-1}(x)$ the boundary of $\{ (a,b)\in U\}\cap \bar{T(x)}$ 
in $\bar{T(x)}$ is finite (with cardinality uniformly bounded)
and non empty because $\bar{T(x)}$ is closed. We have in 
this way a definable function 
$l_{m-1}:X\times Y_{m-1}\into U\cup \{\infty \}$ 
such that $l_{m-1}(x,a)\in \bar{T(x)}$ iff $a\in Y_{m-1}(x)$ and 
$l_{m-1}(x,a)=\infty$ otherwise. 
Similarly, for each
$x\in X$ and $a\in Y_{m-2}(x)$, the definable set
$l_{m-1}(x, \{(a,b)\in Y_{m-1}(x)\})$ has a finite and non empty 
boundary in $\bar{T(x)}$ and we obtain a definable function
$l_{m-2}:X\times Y_{m-2}\into U\cup \{\infty \}$ such that 
$l_{m-2}(x,a)\in \bar{T(x)}$ iff $a\in Y_{m-2}(x)$ and 
$l_{m-2}(x,a)=\infty$ otherwise. Continuing in this 
way, we see that for each $i\in \{1,\dots ,m-1\}$ there is a
definable function $l_i:X\times Y_i\into U\cup \{\infty \}$ 
such that $l_i(x,a)\in \bar{T(x)}$ iff $a\in Y_i(x)$ and 
$l_i(x,a)=\infty$ otherwise. Now, for each $x\in X$
the definable set 
$l_1(x, Y_1(x))$ has finite and non empty boundary in $\bar{T(x)}$ and
so we get a definable choice $l$ for $\bar{T}$ which by
construction is a strong definable choice.

Now let $O$ be the definable neighbourhood of $1$ in $U$ which has
strong definable choice. And consider the definable family
$S:=\{S(x):x\in X\}$ of non empty definable subsets of $O$ where
$S(x):=\{z\in O:l(x)z\in l(x)O\cap T(x)\}$. Note that if $T(x)=T(y)$ 
then $S(x)=S(y)$. Let $s$ be a strong definable choice for $S$. Then 
clearly, $t:=s\cdot l$ is a strong definable choice for $T$. 
\qed

Corollary \ref{cor homo compact} below was also proved in \cite{pst2}
but assuming that $\N$ has definable Skolem functions and using the
theory of $\bigvee$-definable groups. By theorem \ref{thm definable
choice}, the assumption that $\N$ has definable Skolem function is
unnecessary:

\begin{cor}\label{cor homo compact}
Let $A$ be a definably compact definable abelian group. Then 
the following holds. (1) For every definable abelian group $B$, there
is no infinite definable family of definable homomorphisms from $A$
into $B$ or vice-versa. (2) There is no infinite definable family of
definable subgroups of $A$.
\end{cor}

\pf
(1) Let $\gamma :S\times A\into B$ be an infinite definable
family of definable automorphisms of $A$. Then by lemma 2.17 of 
\cite{pst2} there is $\{a_1,\dots ,a_n\}\subseteq A$ such that for 
$s\in S$, $\gamma (s)$ is determined by its values on this finite set.
Therefore, we can identify $S$ with a definable subset of
$A\times \cdots \times A$ ($n$ times). Now the rest of the proof is 
obtained by  adapting the proof of (1) in \cite{pst2}
and using theorem \ref{thm definable choice}. 

(2) The argument in the proof of corollary 5.2 \cite{pst2} reduces
it to case (1). 
\qed

\end{section}

\begin{section}{Definable rings}\label{section definable rings}

In this section we apply 
our result on definable abelian groups to describe
definable rings. We start by recalling some facts about definable rings.

\medskip
Let $U$ be a definable ring. Then by \cite{p} and \cite{opp} $U$ can be 
equipped with a unique definable manifold structure making the ring into a
topological ring, and definable homomorphisms between definable rings
are topological homomorphisms. In fact, it follows from the results in
\cite{pps1}, that if $\N$ is an o-minimal expansion of a real closed 
field then, $U$ equipped with the above unique definable manifold 
structure is a $C^p$ ring for all $p\in \NN$ and definable 
homomorphisms between definable rings are $C^p$ homomorphisms for all
$p\in \NN$.

It follows from the DCC for definable
groups, that $U$ satisfies the descending chain
condition (DCC) on definable left (resp., right and bi-) ideals. Let
$U^0$ be the definable connected component of zero in the additive
group of $U$. Then $U^0$ is the smallest definable ideal of $U$ of
finite index. We say that $U$ is definably connected if $U^0=U$.

Finally we mention the following result from \cite{ps} which we
generalise below.
Let $U$ be an infinite definable associative
ring without zero divisors. Then $U$ 
is a division ring and there is a one-dimensional definable subring 
$I$ of $U$ which is a real closed field such that $U$ is either
$I$, $I(\sqrt{-1})$, or the ring of quaternions over $I$.

\medskip
We now use our main result (theorem \ref{thm the main theorem}), 
the results from \cite{opp}
about rings definable in o-minimal expansions of real closed fields
and Wedderburn theory to prove the following.

\begin{thm}\label{thm definable rings}
Let $U$ be a definable ring. 
Then there is a definable left ideal $V=K\oplus \bigoplus _{j=1}^mW_j
\oplus \bigoplus _{i=1}^nV'_i\oplus
\bigoplus _{i=1}^nV_i$ of $U$ such that 
$K$ is the definably compact, definably connected
 definable left ideal of $U$ of maximal dimension
and for each $j=1,\dots ,m$ (resp., $i=1,\dots ,n$)
there is a semi-bounded o-minimal expansion $\J$$_j$ of a group (resp.,
an o-minimal expansion $\I$$_i$ of a real closed field) definable
in $\N$ such that there is no definable bijection between a distinct
pair among the $J_j$'s and the $I_i$'s, $W_j$ is a direct product of
copies of the additive group of $\J$$_j$ and has zero multiplication,
$V'_i$ is a direct product of
copies of the linearly bounded one-dimensional torsion-free
$\I$$_i$-definable group and has zero multiplication, 
each $V_i$ is an $\I$$_i$-definable ring such that if
$\bar{V_i}:=V_i/ann_{V_i}V_i$ is non-trivial then 
$\bar{V_i}$ is a finitely generated $I_i$-algebra
(and therefore $I_i$-definable) 
and if it is associative then it is $\I$$_i$-definably isomorphic to 
a finitely generated $I_i$-subalgebra
of some $M_{n_i}(I_i)$ and has a nilpotent finitely generated
ideal $Z_i$ such that $\bar{V_i}/Z_i$ is $\I$$_i$-definably isomorphic to
$\bigoplus _{j=1}^{m_i}M_{k_{i,j}}(D_{i,j})$ where for each 
$j=1,\dots m_i$, $D_{i,j}$ is either
$I_i$, $I_i(\sqrt{-1})$, or the ring of quaternions over $I_i$.
Moreover, $U/V$ is a definably compact definable ring. 
\end{thm}   

\pf
If we consider $U$ as an additive definable group and apply theorem  
\ref{thm the main theorem} then $U$ has a definable subgroup $V=K\times
W_1\times \cdots \times W_m\times V'_1\times V_1\times \cdots \times
V'_n\times V_n$ such that $K$ is the definably compact,
definably connected
definable additive subgroup of $U$ of maximal dimension
and for each $j=1,\dots ,m$
(resp., $i=1,\dots ,n$)
there is a semi-bounded o-minimal expansion $\J$$_j$ of a group (resp.,
an o-minimal expansion $\I$$_i$ of a real closed field) definable
in $\N$ such that there is no definable bijection between a distinct
pair among the $J_j$'s and the $I_i$'s, the additive group
$W_j$ is a direct product of
copies of the additive group of $\J$$_j$,
the additive group $V'_i$ is a direct product of
copies of the linearly bounded one-dimensional torsion-free
$\I$$_i$-definable group and each $V_i$ is an $\I$$_i$-definable
additive group. Moreover, any definable additive subgroup of $U/V$ is
definably compact.

It follows easily from this that $V$, $K$, $W_j$'s,
$V'_i$'s and $V_i$'s  are all definable
left ideals of $U$ and $V=K\oplus \bigoplus _{j=1}^mW_j
\oplus \bigoplus _{i=1}^nV'_i\oplus
\bigoplus _{i=1}^nV_i$.  

We now show that each $W_j$ has zero multiplication 
(the proof is the same for each $V'_i$): we have a group homomorphism
$W_j/ann_{W_j}W_j\into End(W_j)$ of additive groups, where $End(W_j)$
is the group of all $\J$$_j$-definable endomorphisms of $W_j$, which
is clearly isomorphic with $M_{n_j}(\Lambda (\J$$_j))$ where $\Lambda
(\J$$_j)$ is the division ring of all $\J$$_j$-definable endomorphisms
of the additive group of $\J$$_j$. By \cite{ms}, $W_j/ann_{W_j}W_j$ 
must be finite, and because $W_j$ is $\J$$_j$-definably connected we have
$W_j=ann_{W_j}W_j$.

By construction of $\I$$_i$, $V_i$ is a $\I$$_i$-definable ring.
Suppose that $\bar{V_i}$ is non-trivial. The fact that each
$\bar{V_i}$ is $\I$$_i$-definably isomorphic with a
finitely generate $I_i$-algebra  
and that if it is associative then it is $\I$$_i$-definably isomorphic to 
a finitely generated $I_i$-subalgebra of some $M_{n_i}(I_i)$
follows from (the proof of) lemma 4.3 in \cite{opp}, and the rest is
just Wedderburn theory (for details see for example the section on
Wedderburn theory in \cite{ab}). 
\qed

\begin{thm}\label{thm compact rings}
A definably compact, definably
connected definable ring has zero multiplication.
\end{thm}

\pf
This is a corollary of the proof of theorem \ref{thm compact G modules}.
\qed

\begin{defn}\label{defn lie ring}
{\em
Recall that a {\it Lie ring} is an additive group $L$ with a bilinear 
product (called bracket) $[x,y]$ such that for all $x,y,z\in L$
(i) $[x,x]=0$ and (ii) $[[x,y],z]+[[y,z],x]+[[z,x],y]=0$ (Jacobi
identity). $L$ is {\it abelian} if for all $x,y\in L$, $[x,y]=0$.
}
\end{defn}

The following facts are proved exactly as
above (using in fact \ref{fact definable lie rings} the Lie ring
analogue of lemma 4.3 in \cite{opp}).

\begin{fact}\label{fact definable lie rings}
{\em
Let $U$ be a definable Lie ring. 
Then there is a definable left ideal $V=K\oplus \bigoplus _{j=1}^mW_j
\oplus \bigoplus _{i=1}^nV'_i\oplus
\bigoplus _{i=1}^nV_i$ of $U$ such that 
$K$ is the definably compact, definably connected
 definable left ideal of $U$ of maximal dimension
and for each $j=1,\dots ,m$ (resp., $i=1,\dots ,n$)
there is a semi-bounded o-minimal expansion $\J$$_j$ of a group (resp.,
an o-minimal expansion $\I$$_i$ of a real closed field) interpretable
in $\N$ such that there is no definable bijection between a distinct
pair among the $J_j$'s and the $I_i$'s, $W_j$ is a direct product of
copies of the additive group of $\J$$_j$ and is an abelian Lie ring,
$V'_i$ is a direct product of
copies of the linearly bounded one-dimensional torsion-free
$\I$$_i$-definable group and is an abelian Lie ring, 
each $V_i$ is an $\I$$_i$-definable ring such that
$\bar{V_i}:=V_i/ann_{V_i}V_i$ is $\I$$_i$-definably isomorphic to 
a finitely generated Lie subalgebra of some $M_{n_i}(I_i)$.
Moreover, $U/V$ is a definably compact definable Lie ring. 
}
\end{fact}

\begin{fact}\label{fact compact lie rings}
{\em
A definably compact, definably connected definable Lie ring is abelian.
}
\end{fact}

\medskip
\textbf{Acknowledgements.} Part of the work presented here is
contained in the authors DPhil Thesis which was financially
supported by JNICT grant PRAXIS XXI/BD/5915/95.
I would like to thank my thesis adviser
Professor Alex Wilkie and Kobi Peterzil for their constant support.
I would also like to thank the EPSRC for current financial support.
\end{section}

\end{document}